\documentclass[12pt, twoside]{amsart}
\usepackage[english]{babel}
\theoremstyle{plain}
\usepackage{mathrsfs}
\usepackage{dsfont}
\usepackage{amssymb}

\makeatletter
\numberwithin{equation}{section}
\numberwithin{figure}{section}
\newtheorem{thm}{\protect\theoremname}[section]
\theoremstyle{plain}

\theoremstyle{definition}
\newtheorem{defn}[thm]{\protect\definitionname}
\theoremstyle{plain}
\newtheorem{prop}[thm]{\protect\propositionname}
\theoremstyle{plain}
\newtheorem{lem}[thm]{\protect\lemmaname}
\theoremstyle{plain}
\newtheorem{cor}[thm]{\protect\corollaryname}

\usepackage{babel}
\addto\captionsamerican{\renewcommand{\corollaryname}{Corollary}}
\addto\captionsamerican{\renewcommand{\definitionname}{Definition}}
\addto\captionsamerican{\renewcommand{\lemmaname}{Lemma}}
\addto\captionsamerican{\renewcommand{\propositionname}{Proposition}}
\addto\captionsamerican{\renewcommand{\theoremname}{Theorem}}
\addto\captionsbritish{\renewcommand{\corollaryname}{Corollary}}
\addto\captionsbritish{\renewcommand{\definitionname}{Definition}}
\addto\captionsbritish{\renewcommand{\lemmaname}{Lemma}}
\addto\captionsbritish{\renewcommand{\propositionname}{Proposition}}
\addto\captionsbritish{\renewcommand{\theoremname}{Theorem}}
\addto\captionsenglish{\renewcommand{\corollaryname}{Corollary}}
\addto\captionsenglish{\renewcommand{\definitionname}{Definition}}
\addto\captionsenglish{\renewcommand{\lemmaname}{Lemma}}
\addto\captionsenglish{\renewcommand{\propositionname}{Proposition}}
\addto\captionsenglish{\renewcommand{\theoremname}{Theorem}}
\providecommand{\corollaryname}{Corollary}
\providecommand{\definitionname}{Definition}
\providecommand{\lemmaname}{Lemma}
\providecommand{\propositionname}{Proposition}
\providecommand{\theoremname}{Theorem}

\addto\captionsbritish{\renewcommand{\definitionname}{Definition}}
\addto\captionsbritish{\renewcommand{\lemmaname}{Lemma}}
\addto\captionsbritish{\renewcommand{\theoremname}{Theorem}}
\addto\captionsenglish{\renewcommand{\definitionname}{Definition}}
\addto\captionsenglish{\renewcommand{\lemmaname}{Lemma}}
\addto\captionsenglish{\renewcommand{\theoremname}{Theorem}}
\providecommand{\definitionname}{Definition}
\providecommand{\lemmaname}{Lemma}
\providecommand{\theoremname}{Theorem}

\newtheoremstyle{boldremark}
    {\dimexpr\topsep/2\relax} 
    {\dimexpr\topsep/2\relax} 
    {}          
    {}          
    {\bfseries} 
    {.}         
    {.5em}      
    {}          

\theoremstyle{boldremark}
\newtheorem{brem} [thm] {Remark} 

\usepackage[colorlinks,pdfpagelabels,pdfstartview = FitH,bookmarksopen
= true,bookmarksnumbered = true,linkcolor = blue,plainpages =
false,hypertexnames = false,citecolor = blue,pagebackref=false,urlcolor=blue]{hyperref}

\allowdisplaybreaks

\addtocounter{page}{1}

\title{Gradient regularity for strongly singular or degenerate elliptic and parabolic equations
\rule{8cm}{0.4pt}\\[0.5em]
Regolarità del gradiente per equazioni ellittiche e paraboliche fortemente singolari o degeneri
}
\author{Pasquale Ambrosio
}
\address{Dipartimento di Matematica, Università di Bologna, Piazza di Porta S. Donato 5, 40126 Bologna, Italy.
}
\email{pasquale.ambrosio@unibo.it
}
\thanks{Bruno Pini Mathematical Analysis Seminar, Vol.  (2025) pp. 2–34   
\\
Dipartimento di Matematica, Universit\`{a} di Bologna
\\
ISSN 2240-2829}

\begin{document}

\begin{abstract}
{We present recent advances in the regularity theory for
weak solutions to some classes of elliptic and parabolic equations
with strongly singular or degenerate structure. The equations under
consideration satisfy standard $p$-growth and $p$-ellipticity conditions
only outside a ball centered at the origin. In the elliptic setting,
we describe Besov and Sobolev regularity results for suitable nonlinear
functions of the gradient of the weak solutions, covering both the
subquadratic ($1<p<2$) and superquadratic ($p\geq2$) regimes. Analogous
results are obtained in the corresponding parabolic framework, where
we address the higher spatial and temporal differentiability of the
solutions under appropriate assumptions on the data.
}

\medskip \noindent {\sc{Sunto.}} 
{Presentiamo alcuni recenti sviluppi nella teoria della regolarità per soluzioni deboli di alcune classi di equazioni ellittiche e paraboliche fortemente singolari o degeneri. Le equazioni considerate soddisfano condizioni standard di crescita e di ellitticità di ordine $p$, tipiche degli operatori di tipo $p$-Laplaciano, ma soltanto all’esterno di una sfera centrata nell’origine. Nel setting ellittico, descriviamo risultati di regolarità di tipo Besov e Sobolev per opportune funzioni non lineari del gradiente delle soluzioni deboli, coprendo sia il caso subquadratico ($1<p<2$) sia quello superquadratico ($p\geq2$). Risultati analoghi sono ottenuti nel corrispondente contesto parabolico, dove si analizza la maggiore regolarità delle soluzioni rispetto alla variabile spaziale e a quella temporale, sotto opportune ipotesi sui termini noti.
}


\end{abstract}
\maketitle
\markboth{PASQUALE AMBROSIO}{GRADIENT REGULARITY FOR SOME CLASSES OF ELLIPTIC AND PARABOLIC EQUATIONS}

\section{Introduction\label{sec:Intro}}

\noindent $\hspace*{1em}$In this note, we present some recent results
proved in \cite{Amb1,Ambr2,AmbGriPas,AmbPass}, concerning the gradient
regularity of weak solutions to certain classes of strongly singular
or degenerate elliptic and parabolic equations.\\
$\hspace*{1em}$As for the elliptic equation in question, this \foreignlanguage{english}{arises
as the Euler-Lagrange equation of an integral functional of the Calculus
of Variations. The energy density of this functional satisfies standard
$p$-growth and $p$-ellipticity conditions with respect to the gradient
variable, but only outside a ball with radius $\lambda>0$ centered
at the origin.} More precisely, we consider the following pair of
variational problems in duality\footnote{The equation $-\,\mathrm{div}\,\sigma=f$ in problem (\ref{eq:P1}) has to be understood as
\[
\int_{\Omega}\langle\sigma,D\varphi\rangle\,dx\,=\int_{\Omega}f\varphi\,dx\,,\,\,\,\,\,\,\mathrm{for}\,\,\mathrm{every}\,\,\varphi\in C^{1}(\overline{\Omega})\,,
\]
so that it incorporates in the weak sense the homogeneous Neumann boundary condition $\langle\sigma,\nu\rangle=0$ on $\partial\Omega$.}: 
\begin{equation}
\inf_{\sigma\,\in\,L^{p'}\left(\Omega,\mathbb{R}^{n}\right)}\left\{ \int_{\Omega}\mathcal{H}(\sigma(x))\,dx:-\,\mathrm{div}\,\sigma=f,\,\,\langle\sigma,\nu\rangle=0\,\,\mathrm{on}\,\,\partial\Omega\right\} \tag{\ensuremath{\mathrm{P1}}}\label{eq:P1}
\end{equation}
and
\begin{equation}
\sup_{u\,\in\,W^{1,p}\left(\Omega\right)}\left\{ \int_{\Omega}u(x)f(x)\,dx-\int_{\Omega}\mathcal{H}^{*}(Du(x))\,dx\right\} ,\tag{\ensuremath{\mathrm{P2}}}\label{eq:P2}
\end{equation}
where $\Omega$ is a bounded connected open subset of $\mathbb{R}^{n}$
($n\geq2$) with Lipschitz boundary, $p\in(1,\infty)$, $p':=p/(p-1)$,
$f$ is a given function defined over $\overline{\Omega}$ with zero
mean (i.e. $\int_{\Omega}f\,dx=0$), $\nu$ denotes the outer normal
versor of $\partial\Omega$ and the function $\mathcal{H}$ is defined
by
\[
\mathcal{H}(\sigma):=\,\frac{1}{p'}\left|\sigma\right|^{p'}+\lambda\left|\sigma\right|,\,\,\,\,\,\,\,\,\sigma\in\mathbb{R}^{n}.
\]
With such a choice, we get
\[
D\mathcal{H}^{*}(z)=\left(\left|z\right|-\lambda\right)_{+}^{p-1}\frac{z}{\left|z\right|}\,,\,\,\,\,\,\,\,\,z\in\mathbb{R}^{n},
\]
where $\mathcal{H}^{*}$ is the Legendre transform of $\mathcal{H}$
and $\left(\,\cdot\,\right)_{+}$ stands for the positive part (see
\cite{BCS}). It is well known that the Euler-Lagrange equation of
the functional in (\ref{eq:P2}) is given by the PDE
\begin{equation}
\begin{cases}
\begin{array}{cc}
-\,\mathrm{div}\,(D\mathcal{H}^{*}(Du))=f & \mathrm{in}\,\,\Omega,\\
\,\,\,\,\,\,\,\langle D\mathcal{H}^{*}(Du),\nu\rangle=0 & \,\,\,\,\mathrm{on}\,\,\partial\Omega,
\end{array}\end{cases}\label{eq:Eul-Lag-1}
\end{equation}
which has to be meant in the distributional sense. If we assume that
$f\in L^{p'}(\Omega)$ and that the infimum in (\ref{eq:P1}) is finite,
then problem (\ref{eq:P1}) consists in minimizing a strictly convex
and coercive functional on $L^{p'}$ subject to a convex and closed
constraint: therefore, a solution $\sigma_{0}$ to (\ref{eq:P1})
exists and must be unique. Moreover, we recall that by standard convex
duality (see for instance \cite{EkTe}), the values $\inf\,$(\ref{eq:P1})
and $\sup\,$(\ref{eq:P2}) coincide and the primal-dual optimality
condition characterizes the minimizer $\sigma_{0}$ of (\ref{eq:P1})
through the equality 
\begin{equation}
\sigma_{0}(x)=D\mathcal{H}^{*}(Du_{0}(x))\,\,\,\,\,\,\,\,\mathrm{for}\,\,\,\mathrm{a.e.}\,\,\,x\in\Omega,\label{eq:opt-cond}
\end{equation}
where $u_{0}\in W^{1,p}(\Omega)$ is a solution of (\ref{eq:P2}). This is equivalent
to the requirement that $u_{0}$ is a weak solution of the Euler-Lagrange
equation (\ref{eq:Eul-Lag-1}), in the sense that
\[
\int_{\Omega}\langle D\mathcal{H}^{*}(Du_{0}(x)),D\varphi(x)\rangle\,dx\,=\int_{\Omega}f(x)\varphi(x)\,dx,\,\,\,\,\,\mathrm{for}\,\,\mathrm{every}\,\,\varphi\in W^{1,p}(\Omega).
\]
Furthermore, since $f$ has zero mean, using the direct methods of
the Calculus of Variations it is not difficult to show that the dual
problem (\ref{eq:P2}) admits at least one solution $u_{0}$ belonging
to
\[
W_{\diamond}^{1,p}(\Omega):=\left\{ u\in W^{1,p}(\Omega):\,\int_{\Omega}u(x)\,dx=0\right\} 
\]
and such that (\ref{eq:opt-cond}) holds, so that $u_{0}$ is a distributional
solution of the strongly singular or degenerate elliptic equation
\begin{equation}
-\,\mathrm{div}\left((\vert Du\vert-\lambda)_{+}^{p-1}\frac{Du}{\vert Du\vert}\right)=f\,\,\,\,\,\,\,\,\mathrm{in}\,\,\,\Omega,\label{eq:elliptic-1}
\end{equation}
under homogeneous Neumann boundary conditions. We also note that,
in general, if one looks at the solutions $u$ of the above equation,
no more than $C^{0,1}$ regularity should be expected for them: indeed,
every $\lambda$-Lipschitz function $u$ is a solution of equation
(\ref{eq:elliptic-1}) with $f=0$. Moreover, when $p\geq2$ we have\smallskip{}
\[
\frac{(\left|Du\right|-\lambda)_{+}^{p-1}}{\left|Du\right|}\left|\xi\right|^{2}\,\leq\,\langle D^{2}\mathcal{H}^{*}(Du)\,\xi,\xi\rangle\,\leq\,(p-1)(\left|Du\right|-\lambda)_{+}^{p-2}\left|\xi\right|^{2}\,\,\,\,\,\,\,\,\mathrm{for\,\,every}\,\,\xi\in\mathbb{R}^{n},
\]
that is, the Hessian matrix $D^{2}\mathcal{H}^{*}(Du)$ degenerates
in the region $\{\vert Du\vert\leq\lambda\}$.\\
$\hspace*{1em}$One of the topics related to the regularity of solutions
to equations as the one in (\ref{eq:elliptic-1}) concerns the study
of their higher differentiability of both integer and fractional order,
and several results are available in this direction; see, for example,
\cite{Amb1,AmbGriPas,Bra,BraCar,BCS,CGHP,CuGiaGioPa,Ru}.
A common aspect of nonlinear elliptic problems with growth rate $p\geq2$
is that the higher differentiability is proved for a nonlinear function
of the gradient of the weak solutions that takes into account the
growth of the structure function of the equation. Indeed, already
for the $p$-Poisson equation (which is obtained from (\ref{eq:elliptic-1})
by setting $\lambda=0$), the higher differentiability is established
for the function\vspace{-0,5mm}
\[
\mathbb{V}_{p}(Du):=\,\vert Du\vert^{\frac{p-2}{2}}Du\,,
\]
as can be seen in many papers, starting from the pioneering one by
Uhlenbeck \cite{Uhl}. In the case of equation (\ref{eq:elliptic-1}),
this phenomenon persists and higher differentiability results hold
true for the function\vspace{-0,5mm}
\begin{equation}
H_{\frac{p}{2}}(Du):=\left(\left|Du\right|-\lambda\right)_{+}^{\frac{p}{2}}\,\frac{Du}{\left|Du\right|}\,.\label{eq:Hp/2}
\end{equation}
In particular, Brasco, Carlier and Santambrogio \cite{BCS} proved
that $H_{\frac{p}{2}}(Du)\in W_{loc}^{1,2}(\Omega,\mathbb{R}^{n})$
if $p\geq2$ and $f\in W^{1,p'}(\Omega)$.\\
$\hspace*{1em}$In Section \ref{sec:1stPaper}, we shall address the
subquadratic (or \textit{singular}) case $1<p<2$, which had been
neglected in the literature, since extra technical difficulties arise
concerning elliptic regularity whenever we are in this case. This
aspect, which also occurs in the classic $p$-harmonic setting, has
been very well explained in \cite[Section 2.6]{BDW}. More precisely,
in the singular case we present four main results established in the
paper \cite{Amb1}.\\
$\hspace*{1em}$In Section \ref{sec:2ndPaper}, we report three main
results proved in the paper \cite{AmbGriPas}. To be more specific,
there we deal with the local $W^{1,2}$-regularity of a \textit{novel}
nonlinear function of the gradient $Du$ of local weak solutions to
(\ref{eq:elliptic-1}), under the following assumptions: \foreignlanguage{american}{\vspace{0.3cm}
}

\noindent $\hspace*{1em}\bullet$$\,\,\,\,f\in B_{p',1,loc}^{\frac{p-2}{p}}(\Omega)$
$\,\,\,\,$if $2<p<\infty$$\,\,$(see Theorem \ref{thm:theo1-2});\vspace{0.3cm}

\noindent $\hspace*{1em}\bullet$ $\,\,\,\,f\in L_{loc}^{{\frac{np}{n(p-1)+2-p}}}(\Omega)$
$\,\,\,\,$if $1<p\leq2$ $\,\,$(see Theorem \ref{thm:theo2-1}).\\

\noindent $\hspace*{1em}$Let us recall again that, for $\lambda=0$,
equation (\ref{eq:elliptic-1}) turns into the $p$-Poisson equation.
For the weak solutions $u\in W^{1,p}(\Omega)$ of this last equation,
Irving and Koch \cite{Irv} proved that\vspace{-2mm} 
\begin{equation}
\mathbb{V}_{p}(Du)\,\in\,W_{loc}^{1,2}(\Omega,\mathbb{R}^{n})\label{eq:Uhlenbeck-1-1}
\end{equation}
if the datum $f$ belongs to the Besov space $B_{p',1}^{\frac{p-2}{p}}(\Omega)$,
with $p>2$ \footnote{Owing to the identity $\vert\mathbb{V}_{p}(Du)\vert=\vert Du\vert^{\frac{p}{2}}$, the local $W^{1,2}$-regularity of $\mathbb{V}_{p}(Du)$ entails an improved local summability of $Du$, which, in turn, provides the foundation for the local Hölder continuity of the weak solution $u$ itself. A pretty simple proof of the Hölder continuity of $u$ is available when $n\geq3$ and $p>\max\,\{2,n-2\}$. In this case, the reasoning is as follows. Since $D\mathbb{V}_{p}(Du)\in L^2_{loc}(\Omega)$, the Sobolev Embedding Theorem ensures that $\mathbb{V}_{p}(Du)\in L^{2n/(n-2)}_{loc}(\Omega)$, that is $Du\in L^{np/(n-2)}_{loc}(\Omega)$. Now observe that
\[
\frac{np}{n-2} >n\,\,\,\,\,\,\mathrm{when}\,\,\,\,\,\,p>\max\,\{2,n-2\}. 
\]
Then, using the Sobolev Embedding Theorem again, we conclude that $u\in C^{0,\alpha}_{loc}(\Omega)$, with $\alpha=1-(n-2)/p$, since it belongs to $W^{1,\frac{np}{n-2}}_{loc}(\Omega)$ by the Poincaré-Wirtinger inequality.\\
$\hspace*{1em}$As for the existence of the weak  second derivatives of $u$ in $L^2_{loc}(\Omega)$, we refer the reader to \cite[Chapter 4, pages 29 and 30]{LindPL}.}. Their assumption on $f$ is essentially sharp, in the
sense that the above result is false if 
\[
f\in B_{p',1}^{s}(\Omega)\,\,\,\,\,\,\,\,\mathrm{with}\,\,\,\,\,s\,<\,\frac{p-2}{p}\,.
\]
Indeed, Brasco and Santambrogio \cite[Section 5]{BraSan} showed with
an explicit example that condition \eqref{eq:Uhlenbeck-1-1} may not
hold if $f$ belongs to a fractional Sobolev space $W_{loc}^{\sigma,p'}(\mathbb{R}^{n})$
with $0<\sigma<(p-2)/p$, which is continuously embedded into $B_{p',1,loc}^{s}(\mathbb{R}^{n})$
whenever $s\in(0,\sigma)$ (see Lemma \ref{lem:emb} below).\\
\foreignlanguage{english}{$\hspace*{1em}$The main results reported
in Section \ref{sec:2ndPaper} are in the spirit of the ones mentioned
above. Indeed, our primary goal in \cite{AmbGriPas} was to find the
assumptions to impose on the datum $f$ \textit{in the scale of local
Besov or Lebesgue spaces} to obtain the $W^{1,2}$-regularity of a
nonlinear function of the gradient of weak solutions to the widely
degenerate or singular equation }(\ref{eq:elliptic-1})\foreignlanguage{english}{.}
\foreignlanguage{american}{We emphasize that, for $\lambda=0$, our
results }give back those contained in \cite{Clop} and \cite[Remark 1.4]{Irv}.\\
\foreignlanguage{english}{$\hspace*{1em}$}A key tool in the proof
of the above-mentioned Theorem \ref{thm:theo1-2} is the duality of
Besov spaces (see Section \ref{subsec:Besov}). In fact, our approach
has been inspired by \cite{BraSan}, where Brasco and Santambrogio
used a duality-based inequality in the setting of fractional Sobolev
spaces, but limiting themselves to the $p$-Poisson equation.\\
\foreignlanguage{english}{$\hspace*{1em}$}In the case $1<p\leq2$,
it is well known that, already for the $p$-Poisson equation, the
higher differentiability of the weak solutions can be achieved without
assuming any differentiability on $f$ (neither of integer nor of
fractional order), but only a suitable degree of integrability. The
sharp assumption on $f$ in the scale of Lebesgue spaces has been
recently found in \cite{Clop}. The aforementioned Theorem \ref{thm:theo2-1}
simply tells us that a result analogous to \cite[Theorem 1.1]{Clop}
holds when dealing with solutions of widely singular equations.\\
$\hspace*{1em}$In Section \ref{sec:3rdPaper}, we move on to the
parabolic version of equation (\ref{eq:elliptic-1}), that is, we
consider the strongly degenerate parabolic PDE
\begin{equation}
u_{t}-\mathrm{div}\left((\vert Du\vert-\lambda)_{+}^{p-1}\frac{Du}{\vert Du\vert}\right)=\widetilde{f}\,\,\,\,\,\,\,\,\mathrm{in}\,\,\,\Omega_{T}=\Omega\times(0,T)\,,\label{eq:1-1}
\end{equation}

\noindent where $p\geq2$, $\Omega$ is a bounded domain in $\mathbb{R}^{n}$
($n\geq2$), $T>0$ and $\lambda$ is a positive constant. The main
feature of the above equation is that the structure function 
\[
H_{p-1}(\xi):=\begin{cases}
\begin{array}{cc}
(\vert\xi\vert-\lambda)_{+}^{p-1}\frac{\xi}{\left|\xi\right|} & \,\,\,\mathrm{if}\,\,\,\xi\in\mathbb{R}^{n}\setminus\{0\},\\
0 & \mathrm{if}\,\,\,\xi=0,\,\,\,\,\,\,\,\,\,\,\,\,\,\,\,
\end{array}\end{cases}
\]
satisfies standard growth and ellipticity conditions for a growth
rate $p\geq2$, but only outside the ball $\{\vert\xi\vert<\lambda\}$.\\
$\hspace*{1em}$In \cite{AmbPass} we succeeded in proving higher
differentiability results for the weak solutions of equation (\ref{eq:1-1}).
These results will be exactly the subject of Section \ref{sec:3rdPaper}.
For further results available in the literature, we refer the reader
to \cite{AMB-frac,AmbBau,BDGP-par,GenPas,Strunk} and the references
therein.\\
$\hspace*{1em}$A motivation for studying equations of the type (\ref{eq:1-1})
can be found in \textit{gas filtration problems taking into account
the initial pressure gradient} (see \cite[Section 1.1]{AmbPass} for
a brief explanation).\\
$\hspace*{1em}$As already pointed out in \cite{Ambr2,BDGP-par}, no more than Lipschitz
regularity can be expected for solutions to (\ref{eq:1-1}). In fact, when $\widetilde{f}=0$, any time-independent
$\lambda$-Lipschitz function solves (\ref{eq:1-1}), and even
more, it is a solution of the associated stationary equation (\ref{eq:elliptic-1}) with $f=0$.\\
$\hspace*{1em}$The first result presented in Section \ref{sec:3rdPaper}
establishes the Sobolev spatial regularity of the vector field defined
in (\ref{eq:Hp/2}), where $Du$ now denotes the \textit{spatial gradient}
of the weak solutions to (\ref{eq:1-1}). This result, in turn, implies
the Sobolev time regularity of the weak solutions, under the assumption
that the datum $\widetilde{f}$ belongs to a suitable Lebesgue-Sobolev
parabolic space. Such results have been obtained by adapting the techniques
for the evolutionary $p$-Laplacian to our more degenerate context.
In fact, for less degenerate parabolic problems, these issues have
been widely investigated, as one can see, for example, in \cite{Duzaar,Giann}
(where $\widetilde{f}=0$) and in \cite{Sche}. Moreover, establishing
the Sobolev regularity of the solutions with respect to time, once
the higher differentiability in space has been obtained, is a quite
usual fact in these problems: see, for instance, \cite{Lind1,Lind2,Lind3}.\\
$\hspace*{1em}$The distinguishing feature of equation (\ref{eq:1-1})
is that its principal part behaves like a \textit{$p$-Laplace operator
only at infinity}. Let us briefly summarize a few previous results
on this topic: the regularity of solutions to parabolic problems with
asymptotic structure of $p$-Laplacian type has been studied in \cite{Isernia},
where a BMO regularity has been proved for solutions to asymptotically
parabolic systems in the case $p=2$ and $\widetilde{f}=0$ (see also
\cite{Kuusi}, where the local Lipschitz continuity of weak solutions
with respect to the spatial variable is established). In addition,
we mention the work \cite{Byun}, where the authors consider nonhomogeneous
parabolic problems involving a discontinuous nonlinearity and an asymptotic
regularity in divergence form of $p$-Laplacian type. There, Byun,
Oh and Wang establish a global Calderón-Zygmund estimate by converting
a given asymptotically regular problem to a suitable regular problem.\\
$\hspace*{1em}$However, it is worth noting that, in Section \ref{sec:3rdPaper},
our assumption on $\widetilde{f}$ is weaker than those considered
in the works mentioned above.\\
$\hspace*{1em}$Finally, in Section \ref{sec:4thPaper} we present
three results proved in the paper \cite{Ambr2}. More precisely, there
we deal with the spatial $W^{1,2}$-regularity of a \textit{novel}
nonlinear function of the spatial gradient of local weak solutions
to equation (\ref{eq:1-1}), under the following assumptions: \foreignlanguage{american}{\vspace{0.3cm}
}

\noindent $\hspace*{1em}\bullet$$\,\,\,\,\widetilde{f}\in L_{loc}^{p'}\left(0,T;B_{p',1,loc}^{\frac{p-2}{p}}(\Omega)\right)$\foreignlanguage{british}{
$\,\,\,\,$if $2<p<\infty$ $\,\,$(see Theorem \ref{thm:theo1-3});}\vspace{0.3cm}

\noindent $\hspace*{1em}\bullet$\foreignlanguage{british}{$\,\,\,\,\widetilde{f}\in L_{loc}^{2}(\Omega_{T})$
$\,\,\,\,$if $p=2$ $\,\,$(see Theorem \ref{thm:nuovo}).}\\

\noindent Actually, the theorems in Section \ref{sec:4thPaper} can
somewhat be viewed as the parabolic counterpart of the elliptic results
presented in Section \ref{sec:2ndPaper}, in the case $p\geq2$.\\
$\hspace*{1em}$Before stating the main assumptions and results, in
Sections \ref{sec:Preliminari} and \ref{sec:Functional} we gather
the preliminary material, including classical notations, essential
definitions, basic properties of the difference quotients of Sobolev
functions, and relevant facts on the function spaces involved in this
note.

\section{Notation and essential definitions\label{sec:Preliminari}}

\noindent $\hspace*{1em}$We start with a list of classical notations.
We shall denote by $C$ or $c$ a general positive constant. Relevant
dependencies on parameters and special constants will be suitably
emphasized using parentheses or subscripts. The norm used on $\mathbb{R}^{k}$,
\foreignlanguage{american}{$k\in\mathbb{N}$}, will be the standard
Euclidean one and it will be denoted by $\left|\,\cdot\,\right|$.
In particular, for vectors $\xi,\eta\in\mathbb{R}^{k}$, we write
$\langle\xi,\eta\rangle$ for the usual inner product and $\left|\xi\right|:=\langle\xi,\xi\rangle^{\frac{1}{2}}$
for the corresponding Euclidean norm.\\
$\hspace*{1em}$For points in space-time, we use the abbreviations
$z=(x,t)$ and $z_{0}=(x_{0},t_{0})$, for spatial variables $x,x_{0}\in\mathbb{R}^{n}$
and times $t,t_{0}\in\mathbb{R}$. We also denote by 
\[
B_{\rho}(x_{0})\,=\,B(x_{0},\rho)\,=\,\left\{ x\in\mathbb{R}^{n}:\left|x-x_{0}\right|<\rho\right\} 
\]
the $n$-dimensional open ball with radius $\rho>0$ and center $x_{0}\in\mathbb{R}^{n}$.
When the center is not important or is clear from the context, we
shall simply write $B_{\rho}\equiv B_{\rho}(x_{0})$. Different balls
in the same context will be assumed to have the same center. Moreover,
we use the notation\vspace{-4mm} 
\[
Q_{\rho}(z_{0}):=\,B_{\rho}(x_{0})\times(t_{0}-\rho^{2},t_{0}),\,\,\,\,\,\,\,\,\,\,z_{0}=(x_{0},t_{0})\in\mathbb{R}^{n}\times\mathbb{R},\,\,\rho>0,
\]
for the backward parabolic cylinder with vertex $(x_{0},t_{0})$ and
width $\rho$. We shall sometimes omit the dependence on the vertex
when all cylinders share the same vertex.\\
$\hspace*{1em}$Now we introduce the auxiliary function $H_{\gamma}:\mathbb{R}^{n}\rightarrow\mathbb{R}^{n}$
defined by 
\[
H_{\gamma}(\xi):=\begin{cases}
\begin{array}{cc}
(\vert\xi\vert-\lambda)_{+}^{\gamma}\,\frac{\xi}{\left|\xi\right|} & \mathrm{if}\,\,\,\xi\neq0,\\
0 & \mathrm{if}\,\,\,\xi=0,
\end{array}\end{cases}
\]

\noindent where $\lambda\geq0$ and $\gamma>0$ are parameters. In
Sections \ref{sec:1stPaper} and \ref{sec:2ndPaper}, we shall deal
with local weak solutions to the elliptic equation (\ref{eq:elliptic-1}).
In this framework, we define a local weak solution to (\ref{eq:elliptic-1})
as follows:
\begin{defn}
\noindent Let $f\in L_{loc}^{1}(\Omega)$. A function $u\in W_{loc}^{1,p}(\Omega)$
is a \textit{local weak solution} of equation (\ref{eq:elliptic-1})
if and only if the condition
\[
\int_{\Omega}\langle H_{p-1}(Du),D\varphi\rangle\,dx\,=\,\int_{\Omega}f\varphi\,dx
\]
is satisfied for all $\varphi\in C_{0}^{\infty}(\Omega)$.
\end{defn}
\noindent $\hspace*{1em}$Similarly, we define a weak solution to
equation (\ref{eq:1-1}) as follows: 
\begin{defn}

\noindent \label{def:soldeb}Let $\widetilde{f}\in L_{loc}^{1}(\Omega_{T})$.
A function 
\[
u\,\in\,C^{0}\left((0,T);L^{2}(\Omega)\right)\cap L^{p}\left(0,T;W^{1,p}(\Omega)\right)
\]
is a \textit{weak solution} of equation (\ref{eq:1-1}) if and only
if the condition
\[
\int_{\Omega_{T}}\left(u\cdot\partial_{t}\varphi-\langle H_{p-1}(Du),D\varphi\rangle\right)\,dz\,=\,-\int_{\Omega_{T}}\widetilde{f}\varphi\,dz
\]
is satisfied for all $\varphi\in C_{0}^{\infty}(\Omega_{T})$.

\end{defn}

\subsection{Difference quotients\label{sec:DiffOpe}}

Here we recall the definition and some well-known
properties of the difference quotients, which can be found, for example,
in \cite{Giu}.

\begin{defn}

\noindent For every vector-valued function $F:\mathbb{R}^{n}\rightarrow\mathbb{R}^{k}$
the \textit{finite difference operator }in the direction $x_{j}$
is defined by
\[
\tau_{j,h}F(x)=F(x+he_{j})-F(x)\,,
\]
where $h\in\mathbb{R}$, $e_{j}$ is the unit vector in the direction
$x_{j}$ and $j\in\{1,\ldots,n\}$.\\
$\hspace*{1em}$The \textit{difference quotient} of $F$ with respect
to $x_{j}$ is defined for $h\in\mathbb{R}\setminus\{0\}$ by 
\[
\Delta_{j,h}F(x)\,=\,\frac{\tau_{j,h}F(x)}{h}\,.
\]

\end{defn}

\noindent When no confusion arises, we shall omit the index $j$ and
simply write $\tau_{h}$ or $\Delta_{h}$ instead of $\tau_{j,h}$
or $\Delta_{j,h}$, respectively. 

\begin{prop}

\noindent Let $\Omega\subset\mathbb{R}^{n}$ be an open set and let
$F\in W^{1,q}(\Omega)$, with $q\geq1$. Moreover, let $G:\Omega\rightarrow\mathbb{R}$
be a measurable function and consider the set
\[
\Omega_{\vert h\vert}:=\left\{ x\in\Omega:\mathrm{dist}\left(x,\partial\Omega\right)>\vert h\vert\right\} .
\]
Then:\\
\\
$\mathrm{(}\mathrm{i}\mathrm{)}$ $\Delta_{h}F\in W^{1,q}(\Omega_{\vert h\vert})$
and $\partial_{x_{i}}(\Delta_{h}F)=\Delta_{h}(\partial_{x_{i}}F)$
for every $\,i\in\{1,\ldots,n\}$.\\

\noindent $\mathrm{(}\mathrm{ii}\mathrm{)}$ If at least one of the
functions $F$ or $G$ has support contained in $\Omega_{\vert h\vert}$,
then
\[
\int_{\Omega}F\,\Delta_{h}G\,dx\,=\,-\int_{\Omega}G\,\Delta_{-h}F\,dx\,.
\]
$\mathrm{(}\mathrm{iii}\mathrm{)}$ We have 
\[
\Delta_{h}(FG)(x)=F(x+he_{j})\Delta_{h}G(x)\,+\,G(x)\Delta_{h}F(x)\,.
\]

\end{prop}

\noindent The next result about the finite difference operator is
a kind of integral version of the Lagrange theorem and can be obtained
by combining \cite[Lemma 8.1]{Giu} with \cite[theorem on page 3]{Maz'ya}.

\begin{lem}

\noindent \label{lem:Giusti1} If $0<\rho<R$, $\vert h\vert<\frac{R-\rho}{2}$,
$1<q<\infty$ and $F\in L_{loc}^{1}(B_{R},\mathbb{R}^{k})$ is such
that $DF\in L^{q}(B_{R},\mathbb{R}^{k\times n})$, then
\[
\int_{B_{\rho}}\left|\tau_{h}F(x)\right|^{q}dx\,\leq\,c^{q}\,\vert h\vert^{q}\int_{B_{R}}\left|DF(x)\right|^{q}dx\,,
\]
where $c$ is a positive constant depending only on $n$. Moreover,
if $F\in L^{q}(B_{R},\mathbb{R}^{k})$, then we have
\[
\int_{B_{\rho}}\left|F(x+he_{j})\right|^{q}dx\,\leq\,\int_{B_{R}}\left|F(x)\right|^{q}dx\,.
\]

\end{lem}

\noindent Finally, we recall the following fundamental result, whose
proof can be found in \cite[Lemma 8.2]{Giu}.

\begin{lem}

\noindent \label{lem:RappIncre} Let $F:\mathbb{R}^{n}\rightarrow\mathbb{R}^{k}$
be a function in $L^{q}(B_{R},\mathbb{R}^{k})$, with $1<q<\infty$.
Assume that there exist $\rho\in(0,R)$ and a constant $M>0$ such
that 
\[
\sum_{j=1}^{n}\int_{B_{\rho}}\left|\tau_{j,h}F(x)\right|^{q}dx\,\leq\,M^{q}\,\vert h\vert^{q}
\]
for all $h\in\mathbb{R}$ satisfying $\vert h\vert<\frac{R-\rho}{2}$.
Then $F\in W^{1,q}(B_{\rho},\mathbb{R}^{k})$ and 
\[
\Vert DF\Vert_{L^{q}(B_{\rho})}\leq M\,.
\]
Moreover, for each $j\in\{1,\ldots,n\}$, 
\[
\Delta_{j,h}F\rightarrow\partial_{x_{j}}F\,\,\,\,\,\,\mathit{in}\,\,L_{loc}^{q}(B_{R},\mathbb{R}^{k})\,\,\,\,\mathit{as}\,\,h\rightarrow0\,.
\]

\end{lem}

\section{Function spaces\label{sec:Functional}}

\noindent $\hspace*{1em}$In this section, we recall the definitions
and basic properties of some function spaces that will be used throughout
this note. We begin with Besov spaces, and then move on to fractional
Sobolev spaces (also known as \textit{Sobolev-Slobodeckij spaces}).

\subsection{Besov spaces\label{subsec:Besov}}

We denote by $\mathcal{S}(\mathbb{R}^{n})$
and $\mathcal{S}'(\mathbb{R}^{n})$ the Schwartz space and the space
of tempered distributions on $\mathbb{R}^{n}$, respectively. If $v\in\mathcal{S}(\mathbb{R}^{n})$,
then 
\begin{equation}
\hat{v}(\xi)=(\mathscr{F}v)(\xi)=(2\pi)^{-n/2}\int_{\mathbb{R}^{n}}e^{-i\,\langle x,\xi\rangle}\,v(x)\,dx,\,\,\,\,\,\,\,\,\xi\in\mathbb{R}^{n},\label{eq:Fourier}
\end{equation}
denotes the Fourier transform of $v$. As usual, $\mathscr{F}^{-1}v$
and $v^{\vee}$ stand for the inverse Fourier transform, given by
the right-hand side of (\ref{eq:Fourier}) with $i$ in place of $-i$.
Both $\mathscr{F}$ and $\mathscr{F}^{-1}$ are extended to $\mathcal{S}'(\mathbb{R}^{n})$
in the standard way.\\
Now, let $\Gamma(\mathbb{R}^{n})$ be the collection of all sequences
$\varphi=\{\varphi_{j}\}_{j=0}^{\infty}\subset\mathcal{S}(\mathbb{R}^{n})$
such that 
\[
\begin{cases}
\begin{array}{cc}
\mathrm{supp}\,\varphi_{0}\subset\{x\in\mathbb{R}^{n}:|x|\le2\}\quad\quad\quad\quad\,\,\,\\
\mathrm{supp}\,\varphi_{j}\subset\{x\in\mathbb{R}^{n}:2^{j-1}\le|x|\le2^{j+1}\} & \mathrm{if}\,\,j\in\mathbb{N},
\end{array}\end{cases}
\]
for every multi-index $\beta$ there exists a positive number $c_{\beta}$
such that 
\[
2^{j\vert\beta\vert}\,\vert D^{\beta}\varphi_{j}(x)\vert\le c_{\beta}\,,\,\,\,\,\,\,\,\,\forall\,j\in\mathbb{N}_{0}\,,\ \forall\,x\in\mathbb{R}^{n}
\]
and 
\[
\sum_{j=0}^{\infty}\varphi_{j}(x)=1\,,\,\,\,\,\,\,\,\,\forall\,x\in\mathbb{R}^{n}.
\]
Then, it is well known that $\Gamma(\mathbb{R}^{n})$ is not empty
(see \cite[Section 2.3.1, Remark 1]{Tri}). Moreover, if $\{\varphi_{j}\}_{j=0}^{\infty}\in\Gamma(\mathbb{R}^{n})$,
the entire analytic functions $(\varphi_{j}\,\hat{v})^{\vee}(x)$
make sense pointwise in $\mathbb{R}^{n}$ for any $v\in\mathcal{S}'(\mathbb{R}^{n})$.
Therefore, the following definition makes sense: 

\begin{defn}

\noindent Let $s\in\mathbb{R}$, $1\leq p,q\le\infty$ and $\varphi=\{\varphi_{j}\}_{j=0}^{\infty}\in\Gamma(\mathbb{R}^{n})$.
We define the \textit{Besov space} $B_{p,q}^{s}(\mathbb{R}^{n})$
as the set of all $v\in\mathcal{S}'(\mathbb{R}^{n})$ such that 
\begin{equation}
\Vert v\Vert_{B_{p,q}^{s}(\mathbb{R}^{n})}:=\left(\sum_{j=0}^{\infty}2^{jsq}\,\Vert(\varphi_{j}\,\hat{v})^{\vee}\Vert_{L^{p}(\mathbb{R}^{n})}^{q}\right)^{\frac{1}{q}}<+\infty\,\,\,\,\,\,\,\,\,\,\mathrm{if}\,\,q<\infty,\label{eq:quasi-norm}
\end{equation}
and 
\begin{equation}
\Vert v\Vert_{B_{p,q}^{s}(\mathbb{R}^{n})}:=\,\sup_{j\,\in\,\mathbb{N}_{0}}\,2^{js}\,\Vert(\varphi_{j}\,\hat{v})^{\vee}\Vert_{L^{p}(\mathbb{R}^{n})}<+\infty\,\,\,\,\,\,\,\,\,\,\mathrm{if}\,\,q=\infty.\label{eq:quasi-norm2}
\end{equation}
\smallskip{}
\begin{brem} The space $B_{p,q}^{s}(\mathbb{R}^{n})$ defined above
is a Banach space with respect to the norm $\Vert\cdot\Vert_{B_{p,q}^{s}(\mathbb{R}^{n})}$.
Obviously, this norm depends on the chosen sequence $\varphi\in\Gamma(\mathbb{R}^{n})$,
but this is not the case for the spaces $B_{p,q}^{s}(\mathbb{R}^{n})$
themselves, in the sense that two different choices for the sequence
$\varphi$ give rise to equivalent norms (see \cite[Sections 2.3.2 and 2.3.3]{Tri}).
This justifies our omission of the dependence on $\varphi$ in the
left-hand side of (\ref{eq:quasi-norm})$-$(\ref{eq:quasi-norm2})
and in the sequel.\end{brem}

\end{defn}

\noindent $\hspace*{1em}$The norms of the \textit{classical Besov
spaces} $B_{p,q}^{s}(\mathbb{R}^{n})$ with $s\in(0,1)$, $1\leq p<\infty$
and $1\le q\le\infty$ can be characterized via differences of the
functions involved (cf. \cite[Section 2.5.12, Theorem 1]{Tri}). More
precisely, for $h\in\mathbb{R}^{n}$ and a measurable function $v:\mathbb{R}^{n}\rightarrow\mathbb{R}$,
let us define 
\[
\delta_{h}v(x):=\,v(x+h)-v(x)\,.
\]
Then we have the equivalence
\[
\Vert v\Vert_{B_{p,q}^{s}(\mathbb{R}^{n})}\,\approx\,\Vert v\Vert_{L^{p}(\mathbb{R}^{n})}\,+\,[v]_{B_{p,q}^{s}(\mathbb{R}^{n})}\,,
\]

\noindent where 
\begin{equation}
[v]_{B_{p,q}^{s}(\mathbb{R}^{n})}:=\biggl({\displaystyle \int_{\mathbb{R}^{n}}\biggl({\displaystyle \int_{\mathbb{R}^{n}}\dfrac{|\delta_{h}v(x)|^{p}}{|h|^{sp}}\,dx\biggr)^{\frac{q}{p}}\dfrac{dh}{|h|^{n}}\biggr)^{\frac{1}{q}},\,\,\,\,\,\,\,\,\,\,\text{if}\,\,\,1\le q<\infty},}\label{eq:BeNorm1}
\end{equation}

\noindent and 
\begin{equation}
[v]_{B_{p,\infty}^{s}(\mathbb{R}^{n})}:={\displaystyle \sup_{h\,\in\,\mathbb{R}^{n}}\biggl({\displaystyle \int_{\mathbb{R}^{n}}\dfrac{|\delta_{h}v(x)|^{p}}{|h|^{sp}}\,dx\biggr)^{\frac{1}{p}}}}.\label{eq:BeNorm2}
\end{equation}

\noindent In (\ref{eq:BeNorm1}), if one simply integrates for $\vert h\vert<r$
for a fixed $r>0$, then an equivalent norm is obtained, since

\noindent 
\[
\biggl({\displaystyle \int_{\{|h|\,\geq\,r\}}\biggl({\displaystyle \int_{\mathbb{R}^{n}}\dfrac{|\delta_{h}v(x)|^{p}}{|h|^{sp}}\,dx\biggr)^{\frac{q}{p}}\dfrac{dh}{|h|^{n}}\biggr)^{\frac{1}{q}}\leq\,c(n,s,p,q,r)\,\Vert v\Vert_{L^{p}(\mathbb{R}^{n})}}}\,.
\]
Similarly, in (\ref{eq:BeNorm2}) one can simply take the supremum
over $|h|\leq r$ and obtain an equivalent norm. By construction,
$B_{p,q}^{s}(\mathbb{R}^{n})\subset L^{p}(\mathbb{R}^{n})$.\medskip{}

\noindent \begin{brem}\label{remBes} For $s\in(0,1)$ and $1\leq p,q<\infty$,
we can simply say that $v\in B_{p,q}^{s}(\mathbb{R}^{n})$ if and
only if $v\in L^{p}(\mathbb{R}^{n})$ and $\frac{\delta_{h}v}{\vert h\vert^{s}}\in L^{q}\left(\frac{dh}{\vert h\vert^{n}};L^{p}(\mathbb{R}^{n})\right)$.\end{brem}\smallskip{}

\noindent $\hspace*{1em}$Let $\Omega$ be an arbitrary open set in
$\mathbb{R}^{n}$. As usual, $\mathcal{D}(\Omega)=C_{0}^{\infty}(\Omega)$
stands for the space of all infinitely differentiable functions in
$\mathbb{R}^{n}$ with compact support in $\Omega$. Let $\mathcal{D}'(\Omega)$
be the dual space of all distributions in $\Omega$ and let $g\in\mathcal{S}'(\mathbb{R}^{n})$.
Then we denote by $g\vert_{\Omega}$ its restriction to $\Omega$,
i.e. 
\[
g\vert_{\Omega}\in\mathcal{D}'(\Omega):\,\,\,\,\,(g\vert_{\Omega})(\phi)=g(\phi)\,\,\,\,\,\,\,\mathrm{for}\,\,\phi\in\mathcal{D}(\Omega).
\]

\begin{defn}

\noindent Let $\Omega$ be an arbitrary domain in $\mathbb{R}^{n}$
with $\Omega\neq\mathbb{R}^{n}$ and let $s\in\mathbb{R}$, $1\leq p\le\infty$
and $1\leq q\le\infty$. Then 
\[
B_{p,q}^{s}(\Omega):=\left\{ v\in\mathcal{D}'(\Omega):\,v=g\vert_{\Omega}\,\,\,\mathrm{for\,\,some}\,\,g\in B_{p,q}^{s}(\mathbb{R}^{n})\right\} 
\]
and 
\[
\Vert v\Vert_{B_{p,q}^{s}(\Omega)}:=\,\inf\,\Vert g\Vert_{B_{p,q}^{s}(\mathbb{R}^{n})}\,,
\]
where the infimum is taken over all $g\in B_{p,q}^{s}(\mathbb{R}^{n})$
such that $g\vert_{\Omega}=v$. 

\end{defn}

\noindent $\hspace*{1em}$For a bounded $C^{\infty}$-domain $\Omega\subset\mathbb{R}^{n}$,
the classical Besov spaces $B_{p,q}^{s}(\Omega)$ with $s\in(0,1)$,
$1\leq p<\infty$ and $1\le q\le\infty$ can be characterized via
differences of the functions involved. More precisely, we have the
following result (see \cite[Section 5.2.2]{Tri2}).

\begin{thm}

\noindent Let $\Omega$ be a bounded $C^{\infty}$-domain in $\mathbb{R}^{n}$.
Let $s\in(0,1)$, $1\leq p<\infty$ and $1\leq q<\infty$. Then, a
function $v:\Omega\rightarrow\mathbb{R}$ belongs to the Besov space
$B_{p,q}^{s}(\Omega)$ if and only if $v\in L^{p}(\Omega)$ and
\begin{equation}
[v]_{B_{p,q}^{s}(\Omega)}:=\left(\int_{\mathbb{R}^{n}}\left(\int_{\Omega}\frac{\left|\delta_{h}v(x)\right|^{p}}{\left|h\right|^{sp}}\cdot\mathds{1}_{\Omega}(x+h)\,dx\right)^{\frac{q}{p}}\frac{dh}{\left|h\right|^{n}}\right)^{\frac{1}{q}}<+\infty.\label{eq:Bes1-1}
\end{equation}
Moreover, the Besov space $B_{p,\infty}^{s}(\Omega)$ consists of
all functions $v\in L^{p}(\Omega)$ such that 
\begin{equation}
[v]_{B_{p,\infty}^{s}(\Omega)}:=\sup_{h\,\in\,\mathbb{R}^{n}}\left(\int_{\Omega}\frac{\left|\delta_{h}v(x)\right|^{p}}{\left|h\right|^{sp}}\cdot\mathds{1}_{\Omega}(x+h)\,dx\right)^{\frac{1}{p}}<+\infty.\label{eq:Bes2-1}
\end{equation}

\end{thm}

\noindent $\hspace*{1em}$For $s\in(0,1)$, $1\leq p<\infty$ and
$1\le q\le\infty$, we also have the equivalence

\noindent 
\[
\Vert v\Vert_{B_{p,q}^{s}(\Omega)}\,\approx\,\Vert v\Vert_{L^{p}(\Omega)}\,+\,[v]_{B_{p,q}^{s}(\Omega)}\,.
\]

\noindent If one replaces $\mathbb{R}^{n}$ in (\ref{eq:Bes1-1})
by a ball $B(0,r)$ for a fixed $r>0$, then an equivalent norm is
obtained. Similarly, in (\ref{eq:Bes2-1}) one can simply take the
supremum over $\left|h\right|\leq r$ and obtain an equivalent norm.\medskip{}

\noindent $\hspace*{1em}$If $s\in\mathbb{R}$, $1\leq p<\infty$
and $1\leq q<\infty$, then $\mathcal{S}(\mathbb{R}^{n})$ is a dense
subset of $B_{p,q}^{s}(\mathbb{R}^{n})$ (cf. \cite[Theorem 2.3.3]{Tri}).
Consequently, in that case, a continuous linear functional on $B_{p,q}^{s}(\mathbb{R}^{n})$
can be interpreted in the usual way as an element of $\mathcal{S}'(\mathbb{R}^{n})$.
More precisely, $g\in\mathcal{S}'(\mathbb{R}^{n})$ belongs to the
dual space $(B_{p,q}^{s}(\mathbb{R}^{n}))'$ of the space $B_{p,q}^{s}(\mathbb{R}^{n})$
if and only if there exists a positive number $c$ such that 
\[
\vert g(\phi)\vert\leq\,c\,\Vert\phi\Vert_{B_{p,q}^{s}(\mathbb{R}^{n})}\,\,\,\,\,\,\,\,\,\,\mathrm{for\,\,all\,\,}\phi\in\mathcal{S}(\mathbb{R}^{n})\,.
\]
We endow $(B_{p,q}^{s}(\mathbb{R}^{n}))'$ with the natural dual norm,
defined by 
\[
\Vert g\Vert_{(B_{p,q}^{s}(\mathbb{R}^{n}))'}=\,\sup\,\left\{ \vert g(\phi)\vert:\phi\in\mathcal{S}(\mathbb{R}^{n})\,\,\,\mathrm{and}\,\,\,\Vert\phi\Vert_{B_{p,q}^{s}(\mathbb{R}^{n})}\leq1\right\} ,\,\,\,\,\,\,\,\,g\in(B_{p,q}^{s}(\mathbb{R}^{n}))'.
\]
Now we recall the following duality formula, which has to be meant
as an isomorphism of normed spaces (see \cite[Theorem 2.11.2]{Tri}). 

\begin{thm}

\noindent \label{thm:duality00} Let $s\in\mathbb{R}$, $1\le p<\infty$
and $1\le q<\infty$. Then 
\[
(B_{p,q}^{s}(\mathbb{R}^{n}))'\,=\,B_{p',q'}^{-s}(\mathbb{R}^{n})\,,
\]
where $p'=\infty$ if $p=1$ (similarly for $q'$). 

\end{thm}

\noindent \begin{brem} The restrictions $p<\infty$ and $q<\infty$
in Theorem \ref{thm:duality00} are natural, since, if either $p=\infty$
or $q=\infty$, then $\mathcal{S}(\mathbb{R}^{n})$ is not dense in
$B_{p,q}^{s}(\mathbb{R}^{n})$, and the density of $\mathcal{S}(\mathbb{R}^{n})$
in $B_{p,q}^{s}(\mathbb{R}^{n})$ is the basis of our interpretation
of the dual space $(B_{p,q}^{s}(\mathbb{R}^{n}))'$.\end{brem}

\begin{defn}

\noindent For $s\in\mathbb{R}$, $1\le p\leq\infty$ and $1\le q\leq\infty$,
we define $\mathring{B}_{p,q}^{s}(\mathbb{R}^{n})$ as the completion
of $\mathcal{S}(\mathbb{R}^{n})$ in $B_{p,q}^{s}(\mathbb{R}^{n})$
with respect to the norm 
\[
v\mapsto\Vert v\Vert_{B_{p,q}^{s}(\mathbb{R}^{n})}\,.
\]
Of course, in the definition above, only the limit cases $\max\,\{p,q\}=\infty$
are of interest. We shall denote by $(\mathring{B}_{p,q}^{s}(\mathbb{R}^{n}))'$
the topological dual of\textit{ }$\mathring{B}_{p,q}^{s}(\mathbb{R}^{n})$,
which is endowed with the natural dual norm 
\[
\Vert g\Vert_{(\mathring{B}_{p,q}^{s}(\mathbb{R}^{n}))'}=\,\sup\,\left\{ \vert g(\phi)\vert:\phi\in\mathcal{S}(\mathbb{R}^{n})\,\,\,\mathrm{and}\,\,\,\Vert\phi\Vert_{B_{p,q}^{s}(\mathbb{R}^{n})}\leq1\right\} ,\,\,\,\,\,\,\,\,g\in(\mathring{B}_{p,q}^{s}(\mathbb{R}^{n}))'.
\]

\end{defn}

\noindent $\hspace*{1em}$The following duality result can be found
in \cite[Section 2.11.2, Remark 2]{Tri} (see also \cite[pages 121 and 122]{Tri0}). 

\begin{thm}

\noindent \label{duality} Let $s\in\mathbb{R}$, $1\le p\le\infty$
and $1\le q\le\infty$. Then 
\[
(\mathring{B}_{p,q}^{s}(\mathbb{R}^{n}))'\,=\,B_{p',q'}^{-s}(\mathbb{R}^{n})\,,
\]
where $p'=1$ if $p=\infty$ (similarly for $q'$).

\end{thm}

\noindent $\hspace*{1em}$One also has the following version of Sobolev
embeddings (a proof can be found in \cite[Proposition 7.12]{Har}).

\begin{lem}

\noindent \label{lem:BesEmbed} Suppose that $0<\alpha<1$.\\
$\mathrm{(}a\mathrm{)}$ If $1<p<\frac{n}{\alpha}$ and $1\leq q\leq p_{\alpha}^{*}:=\frac{np}{n-\alpha p}$,
then there exists a continuous embedding $B_{p,q}^{\alpha}(\mathbb{R}^{n})\hookrightarrow L^{p_{\alpha}^{*}}(\mathbb{R}^{n})$.\\
$\mathrm{(}b\mathrm{)}$ If $p=\frac{n}{\alpha}$ and $1\leq q\leq\infty$,
then there exists a continuous embedding $B_{p,q}^{\alpha}(\mathbb{R}^{n})\hookrightarrow\mathrm{BMO}(\mathbb{R}^{n})$,
where $\mathrm{BMO}$ denotes the space of functions with bounded
mean oscillations \textup{\cite[Chapter 2]{Giu}}.

\end{lem}

\noindent $\hspace*{1em}$We now recall the following inclusions between
Besov spaces (see \cite[Proposition 7.10 and Formula (7.35)]{Har}).

\begin{lem}

\noindent \label{lem:BesInclu-1} Suppose that $0<\beta<\alpha<1$.\\
$\mathrm{(}a\mathrm{)}$ If $1<p\leq\infty$ and $1\leq q\leq r\leq\infty$,
then $B_{p,q}^{\alpha}(\mathbb{R}^{n})\subset B_{p,r}^{\alpha}(\mathbb{R}^{n})$.\\
$\mathrm{(}b\mathrm{)}$ If $1<p\leq\infty$ and $1\leq q,r\leq\infty$,
then $B_{p,q}^{\alpha}(\mathbb{R}^{n})\subset B_{p,r}^{\beta}(\mathbb{R}^{n})$.\\
$\mathrm{(}c\mathrm{)}$ If $1\leq q\leq\infty$, then $B_{\frac{n}{\alpha},q}^{\alpha}(\mathbb{R}^{n})\subset B_{\frac{n}{\beta},q}^{\beta}(\mathbb{R}^{n})$.

\end{lem}

\noindent $\hspace*{1em}$Combining Lemmas \ref{lem:BesEmbed} and
\ref{lem:BesInclu-1}, we obtain the following Sobolev-type embedding
theorem for Besov spaces $B_{p,\infty}^{\alpha}(\mathbb{R}^{n})$
that are excluded from assumptions $\mathrm{(}a\mathrm{)}$ in Lemma
\ref{lem:BesEmbed}.

\begin{lem}

\noindent \label{lem:lemBes} Suppose that $0<\beta<\alpha<1$. If
$1<p<\frac{n}{\alpha}$, then there exists a continuous embedding
$B_{p,\infty}^{\alpha}(\mathbb{R}^{n})\hookrightarrow L^{p_{\beta}^{*}}(\mathbb{R}^{n})$.

\end{lem}

\noindent $\hspace*{1em}$We can also define local Besov spaces as
follows. Given a domain $\Omega\subset\mathbb{R}^{n}$, we say that
a function $v:\Omega\rightarrow\mathbb{R}$ belongs to the \textit{local
Besov space} $B_{p,q,loc}^{s}(\Omega)$ if $\phi\,v\in B_{p,q}^{s}(\mathbb{R}^{n})$
whenever $\phi\in C_{0}^{\infty}(\Omega)$. It is worth noticing that
one can prove suitable local versions of Lemmas \ref{lem:BesEmbed}
and \ref{lem:BesInclu-1}, by using local Besov spaces.\\
$\hspace*{1em}$The following lemma is a consequence of Remark \ref{remBes}
and its proof can be found in \cite[Lemma 7]{BCGOP}.

\begin{lem}

\noindent Let $s\in(0,1)$ and $1\leq p,q<\infty$. A function $v\in L_{loc}^{p}(\Omega)$
belongs to the local Besov space $B_{p,q,loc}^{s}(\Omega)$ if and
only if 
\[
\bigg\Vert\frac{\delta_{h}v}{\left|h\right|^{s}}\bigg\Vert_{L^{q}\left(\frac{dh}{\left|h\right|^{n}}\,;\,L^{p}(B_{r})\right)}<+\infty
\]
for any ball $B_{r}\subset B_{2r}\Subset\Omega$. Here, the measure
$\frac{dh}{\left|h\right|^{n}}$ is restricted to the ball $B(0,r)$
on the $h$-space.

\end{lem}

\noindent $\hspace*{1em}$The next result represents the local counterpart
of Lemma \ref{lem:lemBes}:

\begin{lem}

\noindent \label{lem:emb2-1} On any domain $\Omega\subset\mathbb{R}^{n}$
we have the continuous embedding $B_{p,\infty,loc}^{\alpha}(\Omega)\hookrightarrow L_{loc}^{r}(\Omega)$
for all $r\in\left[1,\,\frac{np}{n-\alpha p}\right)$, provided $\alpha\in(0,1)$
and $1<p<\frac{n}{\alpha}$.

\end{lem}

\noindent We refer to \cite[Sections 30-32]{Tar} for a proof of this
lemma. In fact, the above statement follows by localizing the corresponding
result proved for functions defined on $\mathbb{R}^{n}$ in \cite{Tar},
by simply using a smooth cut-off function.

\noindent $\hspace*{1em}$For the treatment of parabolic equations,
we now give the following definition.

\begin{defn}

\noindent \label{def:BocBes} Let $\Omega$ be a bounded $C^{\infty}$-domain
in $\mathbb{R}^{n}$. Let $s\in(0,1)$, $1\leq p<\infty$ and $1\leq q\leq\infty$.
Then, we say that a map $g\in L^{p}(\Omega\times(t_{0},t_{1}))$ belongs
to the space $L^{p}\left(t_{0},t_{1};B_{p,q}^{s}(\Omega)\right)$
if and only if 
\[
\int_{t_{0}}^{t_{1}}\left[g(\cdot,t)\right]_{B_{p,q}^{s}(\Omega)}^{p}\,dt\,\leq\infty\,,
\]
where $\left[\,\cdot\,\right]_{B_{p,q}^{s}(\Omega)}$ is defined by
(\ref{eq:Bes1-1}) or (\ref{eq:Bes2-1}).

\end{defn}

\noindent In the sequel, we shall use the corresponding local version
of the space defined above, which will be denoted by the subscript
``\textit{loc}''. More precisely, given a bounded domain $\Omega\subset\mathbb{R}^{n}$,
we write $g\in L_{loc}^{p}\left(0,T;B_{p,q,loc}^{s}(\Omega)\right)$
if and only if $g\in L^{p}\left(t_{0},t_{1};B_{p,q}^{s}(\Omega')\right)$
for all domains $\Omega'\times(t_{0},t_{1})\Subset\Omega_{T}:=\Omega\times(0,T)$,
with $\Omega'$ being a (bounded) $C^{\infty}$-domain. We shall also
use the following notation, which is typical of Bochner spaces: 
\[
\Vert g\Vert_{L^{p}(t_{0},\,t_{1}\,;\,B_{p,q}^{s}(\Omega'))}:=\left(\int_{t_{0}}^{t_{1}}\Vert g(\cdot,t)\Vert_{B_{p,q}^{s}(\Omega')}^{p}\,dt\right)^{\frac{1}{p}}.
\]

\subsection{Fractional Sobolev spaces}

Here we recall the definition and some properties
of the fractional Sobolev spaces that appear in the statement of Theorem
\ref{thm:theo1-3} below (for more details, we refer to \cite{DiNezza}).

\noindent $\hspace*{1em}$Let $\Omega$ be a general, possibly nonsmooth,
open set in $\mathbb{R}^{n}$. For any $s\in(0,1)$ and for any $q\in[1,+\infty)$,
we define the fractional Sobolev space $W^{s,q}(\Omega,\mathbb{R}^{k})$
as follows:
\[
W^{s,q}(\Omega,\mathbb{R}^{k}):=\left\{ v\in L^{q}(\Omega,\mathbb{R}^{k}):\frac{\left|v(x)-v(y)\right|}{\left|x-y\right|^{\frac{n}{q}\,+\,s}}\,\in L^{q}(\Omega\times\Omega)\right\} ,
\]
i.e. an intermediate Banach space between $L^{q}(\Omega,\mathbb{R}^{k})$
and $W^{1,q}(\Omega,\mathbb{R}^{k})$, endowed with the norm 
\[
\Vert v\Vert_{W^{s,q}(\Omega)}:=\left(\int_{\Omega}\left|v\right|^{q}dx\,+\,\int_{\Omega}\int_{\Omega}\frac{\left|v(x)-v(y)\right|^{q}}{\left|x-y\right|^{n\,+\,sq}}\,dx\,dy\right)^{\frac{1}{q}},
\]
where the term
\begin{equation}
\left[v\right]_{W^{s,q}(\Omega)}:=\left(\int_{\Omega}\int_{\Omega}\frac{\left|v(x)-v(y)\right|^{q}}{\left|x-y\right|^{n\,+\,sq}}\,dx\,dy\right)^{\frac{1}{q}}\label{eq:Gagliardo}
\end{equation}
is the so-called \textit{Gagliardo seminorm} of $v$.

\noindent \begin{brem} For every $s\in(0,1)$ and every $q\in[1,\infty)$,
we have $B_{q,q}^{s}(\mathbb{R}^{n})=W^{s,q}(\mathbb{R}^{n})$. In
fact, using the change of variable $y=x+h$ in (\ref{eq:Gagliardo})
with $\Omega=\mathbb{R}^{n}$, one gets the seminorm (\ref{eq:BeNorm1})
with $p=q$.\end{brem}\medskip{}

\noindent $\hspace*{1em}$As in the classic case with $s$ being an
integer, the space $W^{\sigma,q}(\Omega)$ is continuously embedded
in $W^{s,q}(\Omega)$ when $s\leq\sigma$, as shown by the next result
(see \cite[Proposition 2.1]{DiNezza}).

\begin{prop}
\noindent Let $\Omega$ be an open set in $\mathbb{R}^{n}$, let $q\in[1,+\infty)$
and let $0<s\leq\sigma<1$. Then there exists a constant $C=C(n,s,q)\geq1$
such that, for any $v\in W^{\sigma,q}(\Omega)$, we have 
\[
\Vert v\Vert_{W^{s,q}(\Omega)}\,\leq\,C\,\Vert v\Vert_{W^{\sigma,q}(\Omega)}\,.
\]
In particular, $W^{\sigma,q}(\Omega)\subseteq W^{s,q}(\Omega)$.

\end{prop}

\noindent $\hspace*{1em}$As is well known when $s\in\mathbb{N}$,
under certain regularity assumptions on the open set $\Omega\subset\mathbb{R}^{n}$,
any function in $W^{s,q}(\Omega)$ can be extended to a function in
$W^{s,q}(\mathbb{R}^{n})$. Extension results are needed to improve
some embedding theorems, in the classic case as well as in the fractional
one. In this regard, we now give the following definition.

\begin{defn}

\noindent For any $s\in(0,1)$ and any $q\in[1,\infty)$, we say that
an open set $\Omega\subseteq\mathbb{R}^{n}$ is an \textit{extension
domain} for $W^{s,q}$ if there exists a positive constant $C=C(n,q,s,\Omega)$
such that for every function $v\in W^{s,q}(\Omega)$ there exists
$\widetilde{v}\in W^{s,q}(\mathbb{R}^{n})$ with $\widetilde{v}=v$
on $\Omega$ and
\[
\Vert\widetilde{v}\Vert_{W^{s,q}(\mathbb{R}^{n})}\,\leq\,C\,\Vert v\Vert_{W^{s,q}(\Omega)}\,.
\]

\end{defn}

\noindent $\hspace*{1em}$In general, an arbitrary open set may not
be an extension domain for $W^{s,q}$. However, the following result
ensures that every open Lipschitz set $\Omega$ with bounded boundary
is an extension domain for $W^{s,q}$ (a proof can be found in \cite[Theorem 5.4]{DiNezza}). 

\begin{thm}

\noindent \label{thm:Lipschitz} Let $q\in[1,+\infty)$, let $s\in(0,1)$
and let $\Omega\subseteq\mathbb{R}^{n}$ be an open set of class $C^{0,1}$
with bounded boundary. Then $W^{s,q}(\Omega)$ is continuously embedded
in $W^{s,q}(\mathbb{R}^{n})$, namely for any $v\in W^{s,q}(\Omega)$
there exists $\widetilde{v}\in W^{s,q}(\mathbb{R}^{n})$ such that
$\widetilde{v}=v$ on $\Omega$ and 
\[
\Vert\widetilde{v}\Vert_{W^{s,q}(\mathbb{R}^{n})}\,\leq\,C\,\Vert v\Vert_{W^{s,q}(\Omega)}
\]
for some positive constant $C=C(n,q,s,\Omega)$.

\end{thm}

\noindent For more information on the problem of characterizing the
class of sets that are extension domains for $W^{s,q}$, we refer
the interested reader to \cite{Zhou}, where an answer to this question
has been given (see also \cite{Jones} and \cite[Chapters 11 and 12]{Leoni}).\\
$\hspace*{1em}$We now recall the following Sobolev-type embedding
theorem, whose proof can be found in \cite[Theorem 6.7]{DiNezza}.

\begin{thm}

\noindent \label{thm:fracemb} Let $s\in(0,1)$ and $q\in[1,+\infty)$
be such that $sq<n$. Let $\Omega\subseteq\mathbb{R}^{n}$ be an extension
domain for $W^{s,q}$. Then there exists a positive constant $C=C(n,q,s,\Omega)$
such that, for any $v\in W^{s,q}(\Omega)$, we have
\[
\Vert v\Vert_{L^{r}(\Omega)}\,\leq\,C\,\Vert v\Vert_{W^{s,q}(\Omega)}
\]
for any $r\in\left[q,q^{*}\right]$; that is, the space $W^{s,q}(\Omega)$
is continuously embedded in $L^{r}(\Omega)$ for any $r\in\left[q,q^{*}\right]$,
where $q^{*}:=nq/(n-sq)$ is the so-called ``fractional critical
exponent''.\\
$\hspace*{1em}$Moreover, if $\Omega$ is bounded, then the space
$W^{s,q}(\Omega)$ is continuously embedded in $L^{r}(\Omega)$ for
any $r\in\left[1,q^{*}\right]$.

\end{thm}

\noindent \begin{brem} In the critical case $r=q^{*}$, the constant
$C$ in Theorem \ref{thm:fracemb} does not depend on $\Omega$ (see
\cite[Remark 6.8]{DiNezza}).\end{brem}\medskip{}

\noindent $\hspace*{1em}$The next embedding result can be obtained
by combining \cite[Section 2.2.2, Remark 3]{Tri} with \cite[Section 2.3.2, Proposition 2(ii)]{Tri}. 

\begin{lem}

\noindent \label{lem:emb} Let $s\in(0,1)$ and $q\ge1$. Then, for
every $\sigma\in(0,1-s)$ we have the continuous embedding $W_{loc}^{s+\sigma,q}(\mathbb{R}^{n})\hookrightarrow B_{q,1,loc}^{s}(\mathbb{R}^{n})$.

\end{lem}

\noindent $\hspace*{1em}$For the treatment of parabolic equations,
we now give the following definition.

\begin{defn}

\noindent Let $q\in[1,+\infty)$ and let $s\in(0,1)$. We say that
a map\\ $g\in L^{q}(\Omega\times(t_{0},t_{1}),\mathbb{R}^{k})$ belongs
to the space $L^{q}\left(t_{0},t_{1};W^{s,q}(\Omega,\mathbb{R}^{k})\right)$
if and only if
\[
\int_{t_{0}}^{t_{1}}\int_{\Omega}\int_{\Omega}\frac{\left|g(x,t)-g(y,t)\right|^{q}}{\left|x-y\right|^{n\,+\,sq}}\,dx\,dy\,dt\,<\infty.
\]
In the sequel, we will use the corresponding local version of the
space defined above, which will be denoted by the subscript ``\textit{loc}''.
More precisely, we write $g\in L_{loc}^{q}\left(0,T;W_{loc}^{s,q}(\Omega,\mathbb{R}^{k})\right)$
if and only if $g\in L^{q}\left(t_{0},t_{1};W^{s,q}(\Omega',\mathbb{R}^{k})\right)$
for all domains $\Omega'\times(t_{0},t_{1})\Subset\Omega_{T}:=\Omega\times(0,T)$.

\end{defn}

\section{Besov estimates for a class of elliptic equations\label{sec:1stPaper}}

\noindent $\hspace*{1em}$In this section, we are concerned with the
regularity properties of the local weak solutions to equation (\ref{eq:elliptic-1}).
All the results that we report here are proved in the paper \cite{Amb1}.\\
$\hspace*{1em}$As already anticipated in the introduction, we are
interested in the higher differentiability of the vector field 
\[
H_{\frac{p}{2}}(Du):=\left(\left|Du\right|-\lambda\right)_{+}^{\frac{p}{2}}\,\frac{Du}{\left|Du\right|}\,,\,\,\,\,\,\,\,\,\,\,\mathrm{for}\,\,p\in(1,2)\,,
\]
where $u\in W_{loc}^{1,p}(\Omega)$ is a local weak solution of $(\ref{eq:elliptic-1})$.
In \cite{Amb1}, we have established a Besov regularity result for
the function $H_{p/2}(Du)$ under suitable assumptions on the right-hand
side $f$ of (\ref{eq:elliptic-1}). As an easy consequence of this
result, we have also deduced the Besov regularity of the solution
$\sigma_{0}$ to (\ref{eq:P1}), under the same hypotheses on $f$.\\
$\hspace*{1em}$More precisely, the main results of this section are
the following theorems. For notation and definitions we refer to Sections
\ref{sec:Preliminari} and \ref{subsec:Besov}.

\begin{thm}[{\foreignlanguage{british}{\textbf{\cite[Theorem 1.1]{Amb1}}}}]

\noindent \foreignlanguage{british}{\label{thm:mainth}Let $n\geq2$,
$p\in(1,2)$ and $\alpha\in(0,1)$. Moreover, let $u\in W_{loc}^{1,p}(\Omega)$
be a local weak solution of equation $(\ref{eq:elliptic-1})$. Then,
the following implications hold:
\[
f\in W_{loc}^{1,p'}(\Omega)\,\,\,\Rightarrow\,\,\,H_{\frac{p}{2}}(Du)\in B_{2,\infty}^{\frac{1}{3-p}}\,\,\,\,\mathit{locally}\,\,\mathit{in}\,\,\Omega,
\]
\[
f\in B_{p',\infty,loc}^{\alpha}(\Omega)\,\,\,\Rightarrow\,\,\,H_{\frac{p}{2}}(Du)\in B_{2,\infty}^{\min\left\{ \frac{\alpha+1}{2},\,\frac{1}{3-p}\right\} }\,\,\,\,\mathit{locally}\,\,\mathit{in}\,\,\Omega.
\]
Furthermore:}\\
\foreignlanguage{british}{}\\
\foreignlanguage{british}{$\mathrm{(}a\mathrm{)}$ if $f\in W_{loc}^{1,p'}(\Omega)$,
then for any ball $B_{R}\Subset\Omega$ the following estimate 
\[
\underset{B_{R/2}}{\int}\left|\tau_{j,h}H_{\frac{p}{2}}(Du)\right|^{2}dx\,\leq C_{1}\left(\Vert Df\Vert_{L^{p'}(B')}\,\Vert Du\Vert_{L^{p}(B')}\left|h\right|^{2}+\,\Vert Du\Vert_{L^{p}(B')}^{p}\left|h\right|^{\frac{2}{3-p}}\right)
\]
holds true for every $j\in\{1,\ldots,n\}$, for every $h\in\mathbb{R}$
such that $\left|h\right|\leq r_{0}<\frac{1}{2}\,\mathrm{dist}\left(B_{R},\partial\Omega\right)$,
for $B'=B_{R}+B(0,r_{0})$ and a positive constant $C_{1}=C_{1}(n,p,R)$.}\\

\noindent $\mathrm{(}b\mathrm{)}$ If, on the other hand, $f\in B_{p',\infty,loc}^{\alpha}(\Omega)$,
then for any ball $B_{R}\Subset\Omega$ the following estimate 
\[
\underset{B_{R/2}}{\int}\left|\tau_{j,h}H_{\frac{p}{2}}(Du)\right|^{2}dx\,\leq C_{2}\left(\Vert f\Vert_{B_{p',\infty}^{\alpha}(B')}\,\Vert Du\Vert_{L^{p}(B')}\left|h\right|^{\alpha+1}+\,\Vert Du\Vert_{L^{p}(B')}^{p}\left|h\right|^{\frac{2}{3-p}}\right)
\]
holds true for every $j\in\{1,\ldots,n\}$, for every $h\in\mathbb{R}$
such that $\left|h\right|\leq r_{0}<\frac{1}{2}\,\mathrm{dist}\left(B_{R},\partial\Omega\right)$,
for $B'=B_{R}+B(0,r_{0})$ and a positive constant $C_{2}=C_{2}(n,p,R)$.

\end{thm}

\begin{thm}[{\foreignlanguage{british}{\textbf{\cite[Theorem 1.2]{Amb1}}}}]

\noindent \foreignlanguage{british}{\label{thm:main2}Let $n\geq2$,
$p\in(1,2)$ and $\alpha\in(0,1)$. Moreover, let $\sigma_{0}\in L^{p'}(\Omega,\mathbb{R}^{n})$
be the solution of $(\ref{eq:P1})$. Then, the following implications
hold:
\[
f\in W_{loc}^{1,p'}\left(\Omega\right)\,\,\,\Rightarrow\,\,\,\sigma_{0}\in B_{p',\infty}^{\frac{2}{p'(3-p)}}\,\,\,\,\mathit{locally}\,\,\mathit{in}\,\,\Omega,
\]
\[
f\in B_{p',\infty,loc}^{\alpha}\left(\Omega\right)\,\,\,\Rightarrow\,\,\,\sigma_{0}\in B_{p',\infty}^{\min\left\{ \frac{\alpha+1}{p'},\,\frac{2}{p'(3-p)}\right\} }\,\,\,\,\mathit{locally}\,\,\mathit{in}\,\,\Omega.
\]
Furthermore:}\\
\foreignlanguage{british}{}\\
\foreignlanguage{british}{$\mathrm{(}a\mathrm{)}$ if $f\in W_{loc}^{1,p'}(\Omega)$,
then for every solution $u$ of problem $(\ref{eq:P2})$ and every
ball $B_{R}\Subset\Omega$, the following estimate 
\[
\underset{B_{R/2}}{\int}\left|\tau_{j,h}\sigma_{0}\right|^{p'}dx\,\leq C_{1}\left(\Vert Df\Vert_{L^{p'}(B')}\,\Vert Du\Vert_{L^{p}(B')}\left|h\right|^{2}+\,\Vert Du\Vert_{L^{p}(B')}^{p}\left|h\right|^{\frac{2}{3-p}}\right)
\]
holds true for every $j\in\{1,\ldots,n\}$, for every $h\in\mathbb{R}$
such that $\left|h\right|\leq r_{0}<\frac{1}{2}\,\mathrm{dist}\left(B_{R},\partial\Omega\right)$,
for $B'=B_{R}+B(0,r_{0})$ and a positive constant $C_{1}=C_{1}(n,p,R)$.}\\

\noindent $\mathrm{(}b\mathrm{)}$ If, on the other hand, $f\in B_{p',\infty,loc}^{\alpha}(\Omega)$,
then for every solution $u$ of problem $(\ref{eq:P2})$ and every
ball $B_{R}\Subset\Omega$, the following estimate 
\[
\underset{B_{R/2}}{\int}\left|\tau_{j,h}\sigma_{0}\right|^{p'}dx\,\leq C_{2}\left(\Vert f\Vert_{B_{p',\infty}^{\alpha}(B')}\,\Vert Du\Vert_{L^{p}(B')}\left|h\right|^{\alpha+1}+\,\Vert Du\Vert_{L^{p}(B')}^{p}\left|h\right|^{\frac{2}{3-p}}\right)
\]
holds true for every $j\in\{1,\ldots,n\}$, for every $h\in\mathbb{R}$
such that $\left|h\right|\leq r_{0}<\frac{1}{2}\,\mathrm{dist}\left(B_{R},\partial\Omega\right)$,
for $B'=B_{R}+B(0,r_{0})$ and a positive constant $C_{2}=C_{2}(n,p,R)$.

\end{thm}

\noindent \begin{brem}If we take $p=2$ and $f\in W_{loc}^{1,2}(\Omega)$
in the statements $\mathrm{(}a\mathrm{)}$ of Theorems \ref{thm:mainth}
and \ref{thm:main2}, then we get back the full Sobolev regularity
results of \cite[Theorem 4.2 and Corollary 4.3]{BCS}, i.e. there
is no discrepancy between our results and the ones contained in \cite{BCS}.\end{brem}\smallskip{}

\noindent $\hspace*{1em}$As a consequence of Theorem \ref{thm:main2}
and Lemma \ref{lem:emb2-1}, we obtain the following higher integrability
result for $\sigma_{0}$.

\begin{cor}[{\foreignlanguage{british}{\textbf{\cite[Corollary 1.4]{Amb1}}}}]

\noindent \foreignlanguage{british}{\label{cor:corol1}Under the
assumptions of Theorem \ref{thm:main2}, we obtain the following implications:
\[
f\in W_{loc}^{1,p'}(\Omega)\,\,\,\Rightarrow\,\,\,\sigma_{0}\in L_{loc}^{r}(\Omega,\mathbb{R}^{n})\,\,\,\,\,\,\mathit{for\,\,all}\,\,r\in\left[1,\frac{n(3-p)\,p'}{n(3-p)-2}\right),
\]
\[
f\in B_{p',\infty,loc}^{\alpha}(\Omega)\,\,\,\Rightarrow\,\,\,\sigma_{0}\in L_{loc}^{s}(\Omega,\mathbb{R}^{n})\,,
\]
where
\[
s=\begin{cases}
\begin{array}{cc}
\mathit{any\,\,value\,\,in}\,\,\left[1,\frac{np'}{n-\alpha-1}\right),\,\, & \mathit{if}\,\,\,0<\alpha<\frac{p-1}{3-p}\,,\\
\\
any\,\,value\,\,in\,\,\left[1,\frac{n(3-p)\,p'}{n(3-p)-2}\right), & \mathit{if}\,\,\,\frac{p-1}{3-p}\leq\alpha<1.
\end{array}\end{cases}
\]
}

\end{cor}

\noindent $\hspace*{1em}$As a direct consequence of the previous
results, we also get a gain of integrability for the gradient of the
local weak solutions of (\ref{eq:elliptic-1}). Indeed, we have the
following result.

\begin{cor}[{\foreignlanguage{british}{\textbf{\cite[Corollary 1.5]{Amb1}}}}]

\noindent \foreignlanguage{british}{\label{cor:corol3}Under the
assumptions of Theorem \ref{thm:mainth}, we obtain the following
implications: 
\[
f\in W_{loc}^{1,p'}(\Omega)\,\,\,\Rightarrow\,\,\,Du\in L_{loc}^{r}(\Omega,\mathbb{R}^{n})\,\,\,\,\,\,\mathit{for\,\,all}\,\,r\in\left[1,\frac{np\,(3-p)}{n(3-p)-2}\right),
\]
\[
f\in B_{p',\infty,loc}^{\alpha}(\Omega)\,\,\,\Rightarrow\,\,\,Du\in L_{loc}^{s}(\Omega,\mathbb{R}^{n})\,,
\]
where 
\[
s=\begin{cases}
\begin{array}{cc}
any\,\,value\,\,in\,\,\left[1,\frac{np}{n-\alpha-1}\right),\,\, & \mathit{if}\,\,\,0<\alpha<\frac{p-1}{3-p}\,,\\
\\
any\,\,value\,\,in\,\,\left[1,\frac{np\,(3-p)}{n(3-p)-2}\right), & \mathit{if}\,\,\,\frac{p-1}{3-p}\leq\alpha<1.
\end{array}\end{cases}
\]
}

\end{cor}
\medskip{}

\noindent $\hspace*{1em}$In \cite{Amb1}, we have also analyzed the
non-singular case $p\geq2$, by proving that, unlike what may occur
in the singular one, the Besov regularity of the datum $f$ translates
into a Besov regularity for the function $H_{\frac{p}{2}}(Du)$ with
no loss in the order of differentiation. More specifically, we have
established the following result.

\begin{thm}[{\foreignlanguage{british}{\textbf{\cite[Theorem 4.1]{Amb1}}}}]

\noindent \foreignlanguage{british}{\label{thm:p>2}Let $n\geq2$,
$p\in[2,\infty)$, $\alpha\in(0,1)$ and $f\in B_{p',\infty,loc}^{\alpha}(\Omega)$.
Moreover, let $u\in W_{loc}^{1,p}(\Omega)$ be a local weak solution
of equation $(\ref{eq:elliptic-1})$. Then 
\[
H_{\frac{p}{2}}(Du)\,\in\,B_{2,\infty,loc}^{\frac{\alpha+1}{2}}(\Omega,\mathbb{R}^{n})\,.
\]
Furthermore, for any ball $B_{R}\Subset\Omega$, the following estimate
\[
\underset{B_{R/2}}{\int}\left|\tau_{j,h}H_{\frac{p}{2}}(Du)\right|^{2}dx\,\leq c\left(\Vert f\Vert_{B_{p',\infty}^{\alpha}(B')}\,\Vert Du\Vert_{L^{p}(B')}\left|h\right|^{\alpha+1}+\,\Vert Du\Vert_{L^{p}(B')}^{p}\left|h\right|^{2}\right)
\]
holds true for every $j\in\{1,\ldots,n\}$, for every $h\in\mathbb{R}$
such that $\left|h\right|\leq r_{0}<\frac{1}{2}\,\mathrm{dist}\left(B_{R},\partial\Omega\right)$,
for $B'=B_{R}+B(0,r_{0})$ and a positive constant $c=c(n,p,R)$.}

\end{thm}

\noindent $\hspace*{1em}$The Besov regularity of $H_{\frac{p}{2}}(Du)$
established by the above theorem allows us to get a gain of integrability
for $Du$. More precisely, we have the following result.

\begin{cor}[{\foreignlanguage{british}{\textbf{\cite[Corollary 4.2]{Amb1}}}}]

\noindent \foreignlanguage{british}{\label{cor:Cor1}Under the assumptions
of Theorem \ref{thm:p>2}, we get 
\[
H_{\frac{p}{2}}(Du)\in L_{loc}^{r}(\Omega,\mathbb{R}^{n})\,\,\,\,\,\mathit{for\,\,all}\,\,r\in\left[1,\frac{2n}{n-\alpha-1}\right)
\]
and
\[
Du\in L_{loc}^{s}(\Omega,\mathbb{R}^{n})\,\,\,\,\,\mathit{for\,\,all}\,\,s\in\left[1,\frac{np}{n-\alpha-1}\right).
\]
}

\end{cor}

\noindent $\hspace*{1em}$Again, if we come back to the variational
problem (\ref{eq:P1}), we obtain the following higher integrability
result for its solution:

\begin{cor}[{\foreignlanguage{british}{\textbf{\cite[Corollary 4.4]{Amb1}}}}]

\noindent \foreignlanguage{british}{\label{cor:Cor2}Let $n\geq2$,
$p\in[2,\infty)$, $\alpha\in(0,1)$ and $f\in B_{p',\infty,loc}^{\alpha}(\Omega)$.
Moreover, let $\sigma_{0}\in L^{p'}(\Omega,\mathbb{R}^{n})$ be the
solution of $(\ref{eq:P1})$. Then 
\[
\sigma_{0}\in L_{loc}^{r}(\Omega,\mathbb{R}^{n})\,\,\,\,\mathit{for\,\,all}\,\,r\in\left[1,\frac{np'}{n-\alpha-1}\right).
\]
}

\end{cor}

\section{Second-order regularity for a class of elliptic equations\label{sec:2ndPaper}}

\noindent $\hspace*{1em}$In this section, we consider again the local
$W^{1,p}$ solutions of (\ref{eq:elliptic-1}). All the results that
we report here are proved in \cite{AmbGriPas}. There, we have established
the local $W^{1,2}$-regularity of a nonlinear function of the gradient
$Du$ of local weak solutions to equation (\ref{eq:elliptic-1}),
by assuming that\foreignlanguage{american}{\vspace{0.3cm}
}

\noindent $\hspace*{1em}\bullet$$\,\,\,\,f\in B_{p',1,loc}^{\frac{p-2}{p}}(\Omega)$\foreignlanguage{british}{
$\,\,\,\,$if $2<p<\infty$ $\,\,$(see Theorem \ref{thm:theo1-2});}\vspace{0.3cm}

\noindent $\hspace*{1em}\bullet$\foreignlanguage{british}{$\,\,\,\,f\in L_{loc}^{{\frac{np}{n(p-1)+2-p}}}(\Omega)$
$\,\,\,\,$if $1<p\leq2$ $\,\,$(see Theorem \ref{thm:theo2-1}).}\\

\noindent The above-mentioned theorems, in turn, imply the local higher
integrability of $Du$ under the same hypotheses on the function $f$
(cf. Corollary \ref{cor:corollario3}).\\
$\hspace*{1em}$The main results of this
section are in the spirit of those contained in \cite{BraSan,Clop,Irv},
which we have already discussed in Section \ref{sec:Intro}. In order to state
these results, we introduce the auxiliary function\vspace{1mm}
\[
\mathcal{G}_{\lambda}(t):=\int_{0}^{t}\frac{\omega^{\frac{p}{2}\,+\,\frac{1}{p-1}}}{(\omega+\lambda)^{1\,+\,\frac{1}{p-1}}}\,d\omega\,\,\,\,\,\,\,\,\,\,\mathrm{for}\,\,t\geq0\,.
\]
Moreover, for $\xi\in\mathbb{R}^{n}$ we define the following vector-valued
function:\vspace{2mm} 
\[
\mathcal{V}_{\lambda}(\xi):=\begin{cases}
\begin{array}{cc}
\mathcal{G}_{\lambda}((\vert\xi\vert-\lambda)_{+})\,\frac{\xi}{\left|\xi\right|} & \,\,\,\mathrm{if}\,\,\,\vert\xi\vert>\lambda,\\
0 & \,\,\,\mathrm{if}\,\,\,\vert\xi\vert\leq\lambda.
\end{array}\end{cases}
\]
Notice that, for $\lambda=0$, we have
\begin{equation}
\mathcal{V}_{0}(\xi)=\,\frac{2}{p}\,\mathbb{V}_{p}(\xi):=\,\frac{2}{p}\,\vert\xi\vert^{\frac{p-2}{2}}\xi\,.\label{eq:V0}
\end{equation}
At this point, our main results read as follows.

\begin{thm}[{\foreignlanguage{british}{\textbf{\cite[Theorem 1.1]{AmbGriPas}}}}]

\noindent \foreignlanguage{british}{\label{thm:theo1-2}Let $n\geq2$,
$p>2$, $\lambda\geq0$ and $f\in B_{p',1,loc}^{\frac{p-2}{p}}(\Omega)$}\foreignlanguage{english}{.}\foreignlanguage{british}{
Moreover, let }\foreignlanguage{english}{$u\in W_{loc}^{1,p}(\Omega)$
be a local weak solution of equation $\mathrm{(\ref{eq:elliptic-1})}$.
Then 
\[
\mathcal{V}_{\lambda}(Du)\,\in\,W_{loc}^{1,2}(\Omega,\mathbb{R}^{n})\,.
\]
Furthermore, for every pair of concentric balls $B_{r}\subset B_{R}\Subset\Omega$
we have 
\begin{align*}
 & \int_{B_{r/4}}\left|D\mathcal{V}_{\lambda}(Du)\right|^{2}dx\\
 & \,\,\,\,\,\,\,\leq\left(C+\,\frac{C}{r^{2}}\right)\left[1+\lambda^{p}+\,\Vert Du\Vert_{L^{p}(B_{R})}^{p}\,+\Vert f\Vert_{L^{p'}(B_{R})}^{p'}\right]+\,C\,\Vert f\Vert_{B_{p',1}^{\frac{p-2}{p}}(B_{R})}^{p'}
\end{align*}
}\foreignlanguage{british}{for a positive constant $C$ depending
only on $n$, $p$ and $R$.}

\end{thm}

\begin{thm}[{\foreignlanguage{british}{\textbf{\cite[Theorem 1.4]{AmbGriPas}}}}]

\noindent \foreignlanguage{british}{\label{thm:theo2-1}Let $n\geq2$,
$1<p\leq2$, $\lambda\geq0$ and $f\in L_{loc}^{\frac{np}{n(p-1)+2-p}}(\Omega)$}\foreignlanguage{english}{.}\foreignlanguage{british}{
Moreover, let}\foreignlanguage{english}{ $u\in W_{loc}^{1,p}(\Omega)$
be a local weak solution of equation $\mathrm{(\ref{eq:elliptic-1})}$.
Then
\[
\mathcal{V}_{\lambda}(Du)\,\in\,W_{loc}^{1,2}(\Omega,\mathbb{R}^{n})\,.
\]
Furthermore, for every pair of concentric balls $B_{r}\subset B_{R}\Subset\Omega$
we have
\begin{align*}
 & \int_{B_{r/4}}\left|D\mathcal{V}_{\lambda}(Du)\right|^{2}dx\\
 & \,\,\,\,\,\,\,\leq\,\frac{C}{r^{2}}\left[1+\lambda^{p}\,+\,\Vert Du\Vert_{L^{p}(B_{R})}^{p}\,+\Vert f\Vert_{L^{\frac{np}{n(p-1)+2-p}}(B_{R})}^{p'}\right]\,+\,C\,\Vert f\Vert_{L^{\frac{np}{n(p-1)+2-p}}(B_{R})}^{p'}
\end{align*}
}\foreignlanguage{british}{for a positive constant $C$ depending
only on $n$, $p$ and $R$}\foreignlanguage{english}{.}

\end{thm}

\noindent \begin{brem}Looking at (\ref{eq:V0}), one can easily understand
that Theorem \ref{thm:theo1-2} extends the result established in
\cite[Remark 1.4]{Irv} to a class of widely degenerate elliptic equations
with standard growth, under a sharp assumption on the order of differentiation
of $f$.\end{brem}

\noindent $\hspace*{1em}$As an easy consequence of the higher differentiability
results in Theorems \ref{thm:theo1-2} and \ref{thm:theo2-1}, since
the gradient of the solution is bounded in the region $\{\vert Du\vert\leq\lambda\}$
and $\mathcal{G}_{\lambda}(t)\approx t^{p/2}$ for large values of
$t$ (see \cite[Lemma 2.8]{AmbGriPas}), we are able to establish
the following higher integrability result for the gradient of local
weak solutions of (\ref{eq:elliptic-1}):\vspace{-0,5mm}

\begin{cor}[{\foreignlanguage{british}{\textbf{\cite[Corollary 1.5]{AmbGriPas}}}}]

\noindent \foreignlanguage{british}{\label{cor:corollario3} Under
the assumptions of Theorem \ref{thm:theo1-2} or Theorem \ref{thm:theo2-1},
we have 
\[
Du\in L_{loc}^{q}(\Omega,\mathbb{R}^{n}),
\]
where\vspace{-1,5mm} 
\begin{align*}
q=\begin{cases}
\begin{array}{cc}
\text{any value in \ensuremath{[1,\infty)}} & \text{if}\,\,\,n=2,\,\,\,\,\,\,\,\,\,\,\,\\
\frac{np}{n-2} & \text{if}\,\,\,n\ge3.\,\,\,\,\,\,\,\,\,\,\,
\end{array}\end{cases}
\end{align*}
}

\end{cor}

\section{Sobolev regularity for a class of parabolic equations: part I\label{sec:3rdPaper}}

\noindent $\hspace*{1em}$In this section, we are interested in the
regularity properties of the weak solutions $u:\Omega_{T}\rightarrow\mathbb{R}$
to the strongly degenerate parabolic equation (\ref{eq:1-1}). All
the results presented here are proved in \cite{AmbPass}. The first
one, given below, can be considered as the parabolic counterpart of
\cite[Theorem 4.2]{BCS}. For notation and definitions we refer to
Section \ref{sec:Preliminari}.

\begin{thm}[{\foreignlanguage{british}{\textbf{\cite[Theorem 1.1]{AmbPass}}}}]

\noindent \foreignlanguage{british}{\label{thm:main4} Let $n\geq2$,
$p\in[2,\infty)$, $\lambda>0$, $\frac{np\,+\,4}{np\,+\,4\,-n}\leq\vartheta<\infty$
and $\widetilde{f}\in L^{\vartheta}\left(0,T;W^{1,\vartheta}(\Omega)\right)$.
Moreover, assume that 
\[
u\in C^{0}\left((0,T);L^{2}(\Omega)\right)\cap L^{p}\left(0,T;W^{1,p}(\Omega)\right)
\]
is a weak solution of equation $\mathrm{(\ref{eq:1-1})}$. Then
\[
H_{\frac{p}{2}}(Du)\,\in\,L_{loc}^{2}\left(0,T;W_{loc}^{1,2}(\Omega,\mathbb{R}^{n})\right).
\]
Furthermore, the estimate\begin{align*}
\underset{Q_{\varrho/2}(z_{0})}{\int}|DH_{\frac{p}{2}}(Du)|^{2}\,dz\,&\leq\,\,c\left(\lambda\,\Vert D\widetilde{f}\Vert_{L^{\vartheta}(Q_{R_{0}})}+\,\Vert D\widetilde{f}\Vert_{L^{\vartheta}(Q_{R_{0}})}^{\frac{np\,+\,4}{np\,+\,2\,-\,n}}\right)\nonumber\\
&\,\,\,\,\,\,\,+\frac{c}{R^{2}}\left(\Vert Du\Vert_{L^{p}(Q_{R_{0}})}^{p}+\,\Vert Du\Vert_{L^{p}(Q_{R_{0}})}^{2}+\lambda^{p}+\lambda^{2}\right)
\end{align*}holds true for any parabolic cylinder $Q_{\varrho}(z_{0})\subset Q_{R}(z_{0})\subset Q_{R_{0}}(z_{0})\Subset\Omega_{T}$
and a positive constant $c$ depending at most on $n$, $p$, $\vartheta$
and $R_{0}$.}

\end{thm}

\noindent $\hspace*{1em}$As anticipated in Section \ref{sec:Intro},
from Theorem \ref{thm:main4} we can deduce that any weak solution
$u$ of (\ref{eq:1-1}) admits a weak time derivative $u_{t}$, which
belongs to the local Lebesgue space $L_{loc}^{\min\,\{\vartheta,\,p'\}}(\Omega_{T})$.
The idea is roughly as follows. Consider equation (\ref{eq:1-1});
since the above theorem tells us that in a certain pointwise sense
the second spatial derivatives of $u$ exist, we may develop the expression
under the divergence symbol. This will give us an expression that
equals $u_{t}$, from which we get the desired summability of the
time derivative. Such an argument must be made more rigorous. Furthermore,
we also need to make explicit \textit{a priori} local estimates. These
are provided in the following theorem.

\begin{thm}[{\foreignlanguage{british}{\textbf{\cite[Theorem 1.2]{AmbPass}}}}]

\noindent \foreignlanguage{british}{\label{thm:timeregularity}Under
the assumptions of Theorem \ref{thm:main4}, the time derivative of
the solution exists in the weak sense and satisfies
\[
\partial_{t}u\,\in\,L_{loc}^{\min\,\{\vartheta,\,p'\}}(\Omega_{T}).
\]
Furthermore, the estimate\begin{align*}
&\left(\underset{Q_{\varrho/2}(z_{0})}{\int}\left|\partial_{t}u\right|^{\min\,\{\vartheta,\,p'\}}\,dz\right)^{\frac{1}{\min\,\{\vartheta,\,p'\}}}\nonumber\\
&\,\,\,\,\,\,\,\leq\,\,c\,\Vert \widetilde{f}\Vert_{L^{\vartheta}(Q_{R_{0}})}+\,c\,\,\Vert Du\Vert_{L^{p}(Q_{R_{0}})}^{\frac{p-2}{2}}\,\left(\lambda\,\Vert D\widetilde{f}\Vert_{L^{\vartheta}(Q_{R_{0}})}+\,\Vert D\widetilde{f}\Vert_{L^{\vartheta}(Q_{R_{0}})}^{\frac{np\,+\,4}{np\,+\,2\,-\,n}}\right)^{\frac{1}{2}}\nonumber\\
&\,\,\,\,\,\,\,\,\,\,\,\,\,\,+\,\frac{c}{R}\left(\Vert Du\Vert_{L^{p}(Q_{R_{0}})}^{2p-2}+\,\Vert Du\Vert_{L^{p}(Q_{R_{0}})}^{p}+(\lambda^{p}+\lambda^{2})\,\Vert Du\Vert_{L^{p}(Q_{R_{0}})}^{p-2}\right)^{\frac{1}{2}}
\end{align*}}

\noindent holds true for any parabolic cylinder $Q_{\varrho}(z_{0})\subset Q_{R}(z_{0})\subset Q_{R_{0}}(z_{0})\Subset\Omega_{T}$
and a positive constant $c$ depending on $n$, $p$, $\vartheta$
and $R_{0}$. 

\end{thm}

\noindent \begin{brem}It is worth pointing out that, starting from
the weaker assumption 
\begin{equation}
\widetilde{f}\in L^{\frac{np\,+\,4}{np\,+\,4\,-n}}\left(0,T;W^{1,\frac{np\,+\,4}{np\,+\,4\,-n}}(\Omega)\right),\label{eq:f-assu}
\end{equation}
Sobolev regularity results such as those of Theorems \ref{thm:main4}
and \ref{thm:timeregularity} seemed not to have been established
yet for weak solutions to parabolic PDEs that are far less degenerate
than equation (\ref{eq:1-1}). In particular, the results of Theorems
\ref{thm:main4} and \ref{thm:timeregularity} can be easily extended
to the case $\lambda=0$, i.e. to the evolutionary $p$-Poisson equation
\begin{equation}
u_{t}-\mathrm{div}\left(\vert Du\vert^{p-2}Du\right)=\widetilde{f}\,\,\,\,\,\,\,\,\,\mathrm{in}\,\,\,\Omega_{T},\label{eq:p-Poisson}
\end{equation}
under the assumption (\ref{eq:f-assu}). Therefore, our results have
improved the existing literature, already for equations of the form
(\ref{eq:p-Poisson}), which exhibit a milder degeneracy.\end{brem}

\section{Sobolev regularity for a class of parabolic equations: part II\label{sec:4thPaper}}

\noindent $\hspace*{1em}$In this final section, we are interested
in the local weak solutions $u:\Omega_{T}\rightarrow\mathbb{R}$ of
the strongly degenerate parabolic equation (\ref{eq:1-1}). All the
results presented here are proved in \cite{Ambr2}. There, we have
established the spatial $W^{1,2}$-regularity of a nonlinear function
of the spatial gradient $Du$ of the local weak solutions to (\ref{eq:1-1}),
by assuming that \foreignlanguage{american}{\vspace{0.3cm}
}

\noindent $\hspace*{1em}\bullet$$\,\,\,\,\widetilde{f}\in L_{loc}^{p'}\left(0,T;B_{p',1,loc}^{\frac{p-2}{p}}(\Omega)\right)$\foreignlanguage{british}{
$\,\,\,\,$if $p>2$ $\,\,$(see Theorem \ref{thm:theo1-3});}\vspace{0.3cm}

\noindent $\hspace*{1em}\bullet$\foreignlanguage{british}{$\,\,\,\,\widetilde{f}\in L_{loc}^{2}(\Omega_{T})$
$\,\,\,\,$if $p=2$ $\,\,$(see Theorem \ref{thm:nuovo}).}\\

\noindent The above-mentioned theorems, in turn, imply the Sobolev
time regularity of the local weak solutions to the evolutionary $p$-Poisson
equation, under the same hypotheses on the function $\widetilde{f}$
(cf. Theorem \ref{thm:theo2-2}, where we only address the case $p>2$,
since the Sobolev time regularity is well known for the heat equation
with source term in $L_{loc}^{2}(\Omega_{T})$).\\
$\hspace*{1em}$As already noted in Section \ref{sec:Intro}, the first two
theorems of this section can somewhat be considered as the parabolic
analog of the elliptic results presented in Section \ref{sec:2ndPaper},
in the case $p\geq2$. In order to state these theorems, we introduce
the auxiliary function
\begin{equation}
\mathcal{G}_{\alpha,\lambda}(t):=\int_{0}^{t}\frac{\omega^{\frac{p-1+2\alpha}{2}}}{(\omega+\lambda)^{\frac{1+2\alpha}{2}}}\,d\omega\,\,\,\,\,\,\,\,\,\,\mathrm{for}\,\,t\geq0\,,\label{eq:Gfun}
\end{equation}
where $\alpha\geq0$. Moreover, for $\xi\in\mathbb{R}^{n}$ we define
the following vector field: 
\begin{equation}
\mathcal{V}_{\alpha,\lambda}(\xi):=\begin{cases}
\begin{array}{cc}
\mathcal{G}_{\alpha,\lambda}((\vert\xi\vert-\lambda)_{+})\,\frac{\xi}{\left|\xi\right|} & \,\,\,\mathrm{if}\,\,\,\vert\xi\vert>\lambda,\\
0 & \,\,\,\mathrm{if}\,\,\,\vert\xi\vert\leq\lambda.
\end{array}\end{cases}\label{eq:Vfun}
\end{equation}
At this point, our main results read as follows.

\begin{thm}[{\foreignlanguage{british}{\textbf{\cite[Theorem 1.1]{Ambr2}}}}]

\noindent \foreignlanguage{british}{\label{thm:theo1-3} Let $n\geq2$,
$p>2$, $\lambda\geq0$ and $\widetilde{f}\in L_{loc}^{p'}\left(0,T;B_{p',1,loc}^{\frac{p-2}{p}}(\Omega)\right)$.
Moreover, let 
\begin{equation}
\alpha=\begin{cases}
\begin{array}{cc}
0 & \textit{if}\,\,\,\lambda=0,\\
\textit{any\,\,value\,\,in}\,\,\left[\frac{p+1}{2(p-1)},\infty\right) & \textit{if}\,\,\,\lambda>0,
\end{array}\end{cases}\label{eq:alpha-1}
\end{equation}
and assume that 
\[
u\in C^{0}\left((0,T);L^{2}(\Omega)\right)\cap L_{loc}^{p}\left(0,T;W_{loc}^{1,p}(\Omega)\right)
\]
is a local weak solution of equation $\mathrm{(\ref{eq:1-1})}$. Then
\[
\mathcal{V}_{\alpha,\lambda}(Du)\,\in\,L_{loc}^{2}\left(0,T;W_{loc}^{1,2}(\Omega,\mathbb{R}^{n})\right),
\]
where the function $\mathcal{V}_{\alpha,\lambda}$ is defined according
to $(\ref{eq:Gfun})$$-$$(\ref{eq:Vfun})$. Furthermore, for any
parabolic cylinder $Q_{r}(z_{0})\subset Q_{\rho}(z_{0})\subset Q_{R}(z_{0})\Subset\Omega_{T}$
we have\begin{align*}
\int_{Q_{r/2}(z_{0})}\left|D_{x}\mathcal{V}_{\alpha,\lambda}(Du)\right|^{2}dz\,&\leq \left(C\,+\,\frac{C}{\rho^{2}}\right)\left[\Vert Du\Vert_{L^{p}(Q_{R})}^{p}+\Vert Du\Vert_{L^{p}(Q_{R})}^{2}+\lambda^{p}+\lambda^{2}+1\right]\nonumber\\
&\,\,\,\,\,\,\,+\,C\,\Vert \widetilde{f}\Vert_{L^{p'}\left(t_{0}-R^{2},t_{0};B_{p',1}^{\frac{p-2}{p}}(B_{R}(x_{0}))\right)}^{p'}
\end{align*} for a positive constant $C$ depending only on $n$, $p$ and $R$
in the case $\lambda=0$, and additionally on $\alpha$ if $\lambda>0$.
Besides, if $\lambda=0$ we get 
\[
Du\,\in\,L_{loc}^{p}\left(0,T;W_{loc}^{\sigma,p}(\Omega,\mathbb{R}^{n})\right)\,\,\,\,\,\,\,\,\,\mathit{for}\,\,\mathit{all}\,\,\sigma\in\left(0,\frac{2}{p}\right).
\]
}

\end{thm}

\begin{thm}[{\foreignlanguage{british}{\textbf{\cite[Theorem 1.3]{Ambr2}}}}]

\noindent \foreignlanguage{british}{\label{thm:nuovo} Let $n\geq2$,
$\lambda\geq0$ and $\widetilde{f}\in L_{loc}^{2}(\Omega_{T})$. Moreover,
let 
\begin{equation}
\alpha=\begin{cases}
\begin{array}{cc}
0 & \textit{if}\,\,\,\lambda=0,\\
\textit{any\,\,value\,\,in}\,\,\left[\frac{3}{2},\infty\right) & \textit{if}\,\,\,\lambda>0,
\end{array}\end{cases}\label{eq:alpha-2}
\end{equation}
and assume that
\[
u\in C^{0}\left((0,T);L^{2}(\Omega)\right)\cap L_{loc}^{2}\left(0,T;W_{loc}^{1,2}(\Omega)\right)
\]
is a local weak solution of the equation 
\[
u_{t}-\mathrm{div}\left((\vert Du\vert-\lambda)_{+}\,\frac{Du}{\vert Du\vert}\right)=\widetilde{f}\,\,\,\,\,\,\,\,\mathit{in}\,\,\,\Omega_{T}.
\]
Then 
\[
\mathcal{V}_{\alpha,\lambda}(Du)\,\in\,L_{loc}^{2}\left(0,T;W_{loc}^{1,2}(\Omega,\mathbb{R}^{n})\right).
\]
Furthermore, for any parabolic cylinder $Q_{r}(z_{0})\subset Q_{\rho}(z_{0})\subset Q_{R}(z_{0})\Subset\Omega_{T}$
we have
\[
\int_{Q_{r/2}(z_{0})}\left|D_{x}\mathcal{V}_{\alpha,\lambda}(Du)\right|^{2}dz\,\leq\,\frac{C}{\rho^{2}}\left(\Vert Du\Vert_{L^{2}(Q_{R})}^{2}+\lambda^{2}+1\right)\,+\,C\left(\Vert\widetilde{f}\Vert_{L^{2}(Q_{R})}^{2}+1\right)
\]
for a positive constant $C$ depending only on $n$ and $R$ in the
case $\lambda=0$, and additionally on $\alpha$ if $\lambda>0$.}

\end{thm}

\noindent \begin{brem}Notice that, for every $\alpha\geq0$, we have
\begin{equation}
\mathcal{V}_{\alpha,0}(\xi)=\,\frac{2}{p}\,\mathbb{V}_{p}(\xi):=\,\frac{2}{p}\,\vert\xi\vert^{\frac{p-2}{2}}\xi\,,\label{eq:V0-1}
\end{equation}
i.e. $\mathcal{V}_{\alpha,0}$ is actually independent of the parameter
$\alpha$, which explains the choices $(\ref{eq:alpha-1})_{1}$ and
$(\ref{eq:alpha-2})_{1}$ in the statements of the above theorems.
The conditions $(\ref{eq:alpha-1})_{2}$ and $(\ref{eq:alpha-2})_{2}$
are instead needed to carry out the proofs of Theorems \ref{thm:theo1-3}
and \ref{thm:nuovo} (see \cite{Ambr2} for more details). Moreover,
looking at (\ref{eq:V0-1}), one can easily understand that, on the
one hand, Theorem \ref{thm:theo1-3} extends the result established
in \cite[Lemma 5.1]{Duzaar} for the parabolic $p$-Laplace equation
to a widely degenerate parabolic setting. On the other hand, it extends
the aforementioned result to the case of data in a suitable Lebesgue-Besov
parabolic space, which turns out to be optimal, as can be seen by
appropriately adapting the example in \cite[Section 5]{BraSan} to
the parabolic context (in this regard, see also \cite[page 3]{AmbGriPas}).\end{brem}

\noindent $\hspace*{1em}$For $p>2$, we now consider the evolutionary
$p$-Poisson equation
\begin{equation}
u_{t}-\mathrm{div}\,(\vert Du\vert^{p-2}Du)=\widetilde{f}\,\,\,\,\,\,\,\,\mathrm{in}\,\,\,\Omega_{T}\,,\label{eq:p-Poisson-1}
\end{equation}
which is obtained from equation (\ref{eq:1-1}) by setting $\lambda=0$.
From Theorem \ref{thm:theo1-3}, one can easily deduce that the local
weak solutions of (\ref{eq:p-Poisson-1}) admit a weak time derivative
which belongs to the local Lebesgue space $L_{loc}^{p'}(\Omega_{T})$.
The idea is essentially the same as in the proof of Theorem \ref{thm:timeregularity}
and leads us to the following result.

\begin{thm}[{\foreignlanguage{british}{\textbf{\cite[Theorem 1.5]{Ambr2}}}}]

\noindent \foreignlanguage{british}{\label{thm:theo2-2}Let $n\geq2$,
$p>2$ and $\widetilde{f}\in L_{loc}^{p'}\left(0,T;B_{p',1,loc}^{\frac{p-2}{p}}(\Omega)\right)$.
Moreover, assume that 
\[
u\in C^{0}\left((0,T);L^{2}(\Omega)\right)\cap L_{loc}^{p}\left(0,T;W_{loc}^{1,p}(\Omega)\right)
\]
is a local weak solution of equation $\mathrm{(\ref{eq:p-Poisson-1})}$.
Then, the time derivative of $u$ exists in the weak sense and satisfies
\[
\partial_{t}u\,\in\,L_{loc}^{p'}(\Omega_{T}).
\]
Furthermore, for any parabolic cylinder $Q_{r}(z_{0})\subset Q_{\rho}(z_{0})\subset Q_{R}(z_{0})\Subset\Omega_{T}$
we have \begin{align*}
\left(\int_{Q_{r/2}(z_{0})}\left|\partial_{t}u\right|^{p'}dz\right)^{\frac{1}{p'}}& \leq \left(C\,+\,\frac{C}{\rho}\right)\left[\Vert Du\Vert_{L^{p}(Q_{R})}^{p-1}+\Vert Du\Vert_{L^{p}(Q_{R})}^{\frac{p}{2}}+\Vert Du\Vert_{L^{p}(Q_{R})}^{\frac{p-2}{2}}\right]\nonumber\\
&\,\,\,\,\,\,\,+\,C\,\Vert Du\Vert_{L^{p}(Q_{R})}^{\frac{p-2}{2}}\,\Vert \widetilde{f}\Vert_{L^{p'}\left(t_{0}-R^{2},t_{0};B_{p',1}^{\frac{p-2}{p}}(B_{R}(x_{0}))\right)}^{\frac{p'}{2}}\,+\,\Vert \widetilde{f}\Vert_{L^{p'}(Q_{R})}
\end{align*}for a positive constant $C$ depending only on $n$, $p$ and $R$.}

\end{thm}

\noindent \begin{brem}It is worth pointing out that, starting from
the weaker assumption
\[
\widetilde{f}\in L_{loc}^{p'}\left(0,T;B_{p',1,loc}^{\frac{p-2}{p}}(\Omega)\right)\,\,\,\,\,\,\,\,\mathrm{with}\,\,\,p>2,
\]
Sobolev regularity results such as those of Theorems \ref{thm:theo1-3}
and \ref{thm:theo2-2} seemed not to have been established yet for
weak solutions to parabolic PDEs that are far less degenerate than
equation (\ref{eq:1-1}) with $\lambda>0$. In particular, the results
of Theorems \ref{thm:theo1-3} and \ref{thm:theo2-2} permit to improve
the existing literature, already for the evolutionary $p$-Poisson
equation (\ref{eq:p-Poisson-1}), which exhibits a milder degeneracy.
Moreover, as far as we know, for $\lambda>0$ the result of Theorem
\ref{thm:nuovo} is completely new.\end{brem}

\noindent \medskip{}

\noindent \textbf{Acknowledgments.} The present note is based on a
talk given by the author at the University of Bologna in October 2025,
in the series \textit{Seminari di Analisi Matematica Bruno Pini}.
The author is a member of the Gruppo Nazionale per l'Analisi Matematica,
la Probabilità e le loro Applicazioni (GNAMPA) of the Istituto Nazionale
di Alta Matematica (INdAM), and is partially supported by the INdAM--GNAMPA
2025 Project ``Regolarità ed esistenza per operatori anisotropi”
(CUP E5324001950001). The author also acknowledges financial
support from the IADE\_CITTI\_2020 Project ``Intersectorial applications
of differential equations'' (CUP J34I20000980006).

\bibliographystyle{alpha}

\begin{thebibliography}{111}


\bibitem{Amb1}P. Ambrosio, \textit{Besov regularity for a class of
singular or degenerate elliptic equations}, J. Math. Anal. Appl. \textbf{505}(2),
125636 (2022).

\bibitem{AMB-frac}P. Ambrosio,\textit{ Fractional Sobolev regularity
for solutions to a strongly degenerate parabolic equation}, Forum
Math. \textbf{35}(6), 1485-1497 (2023).

\bibitem{Ambr2}P. Ambrosio, \textit{Sharp Sobolev regularity for
widely degenerate parabolic equations}, Calc. Var. \textbf{64}, 32
(2025).

\bibitem{AmbBau}P. Ambrosio, F. Bäuerlein, \textit{Gradient bounds
for strongly singular or degenerate parabolic systems}, J. Differ.
Equ.\textbf{ 401}, 492-549 (2024).

\bibitem{AmbGriPas}P. Ambrosio, A.G. Grimaldi, A. Passarelli di Napoli,
\textit{On the second-order regularity of solutions to widely singular
or degenerate elliptic equations}, Annali di Matematica Pura ed Applicata
(2025). DOI: \url{https://doi.org/10.1007/s10231-025-01607-7}.

\bibitem{AmbPass}P. Ambrosio, A. Passarelli di Napoli, \textit{Regularity
results for a class of widely degenerate parabolic equations}, Adv.
Calc. Var. \textbf{17}(3), 805-829 (2023). 

\bibitem{BCGOP}\foreignlanguage{english}{A.L. Baisón, A. Clop, R.
Giova, J. Orobitg, A. Passarelli di Napoli, \textit{Fractional differentiability
for solutions of nonlinear elliptic equations}, Potential Anal.\textbf{
46} (3), 403-430 (2017).}

\bibitem{BDW}\foreignlanguage{english}{A.K. Balci, L. Diening, M.
Weimar, \textit{Higher order Calderón-Zygmund estimates for the $p$-Laplace
equation}, J. Differ. Equ. \textbf{268}, 590-635 (2020).}

\bibitem{BDGP-par}\foreignlanguage{english}{V. Bögelein, F. Duzaar,
R. Giova, A. Passarelli di Napoli, \textit{Gradient regularity for
a class of widely degenerate parabolic systems}, SIAM J. Math. Anal.
\textbf{56}(4), 5017-5078 (2024).}

\bibitem{Bra}\foreignlanguage{english}{L. Brasco, \textit{Global
$L^{\infty}$ gradient estimates for solutions to a certain degenerate
elliptic equation}, Nonlinear Anal. \textbf{74}, 516-531 (2011).}


\bibitem{BraCar}L. Brasco, G. Carlier, \textit{Congested traffic
equilibria and degenerate anisotropic PDEs}, Dyn. Games Appl. \textbf{3}
(4), 508-522 (2013).


\bibitem{BCS}L. Brasco, G. Carlier, F. Santambrogio, \textit{Congested
traffic dynamics, weak flows and very degenerate elliptic equations}
{[}corrected version of mr2584740{]}, J. Math. Pures Appl. (9) \textbf{93}(6),
652-671 (2010).

\bibitem{BraSan}L. Brasco, F. Santambrogio, \textit{A sharp estimate
à la Calderón-Zygmund for the $p$-Laplacian}, Commun. Contemp. Math.
\textbf{20}(3), 1750030 (2018).

\bibitem{Byun}\foreignlanguage{english}{S.-S. Byun, J. Oh, L. Wang,
\textit{Global Calderón-Zygmund theory for asymptotically regular
nonlinear elliptic and parabolic equations}, Int. Math. Res. Not.
IMRN\textbf{ 2015} (2015), no. 17, 8289-8308.}

\bibitem{Clop}\foreignlanguage{english}{A. Clop, A. Gentile, A. Passarelli
di Napoli, \textit{Higher differentiability results for solutions
to a class of non-homogeneous elliptic problems under sub-quadratic
growth conditions}, Bull. Math. Sci. \textbf{13}(12), 2350008 (2023).}

\bibitem{CGHP}\foreignlanguage{english}{A. Clop, R. Giova, F. Hatami,
A. Passarelli di Napoli, \textit{Very degenerate elliptic equations
under almost critical Sobolev regularity}, Forum Math. \textbf{32}
(6), 1515-1537 (2020).}

\bibitem{CuGiaGioPa}\foreignlanguage{english}{G. Cupini, F. Giannetti,
R. Giova, A. Passarelli di Napoli, \textit{Regularity results for
vectorial minimizers of a class of degenerate convex integrals}, J.
Differ. Equ. \textbf{265}, 4375-4416 (2018).}

\bibitem{DiNezza}E. Di Nezza, G. Palatucci, E. Valdinoci, \textit{Hitchhiker's
guide to the fractional Sobolev spaces}, Bull. Sci. Math. \textbf{136},
no. 5, 521-573 (2012).

\bibitem{Duzaar}\foreignlanguage{english}{F. Duzaar, G. Mingione,
K. Steffen, \textit{Parabolic systems with polynomial growth and regularity},
Mem. Amer. Math. Soc. 214 (2011).}

\bibitem{EkTe}\foreignlanguage{english}{I. Ekeland, R. Témam, \textit{Convex
Analysis and Variational Problems}, Classics in Applied Mathematics,
vol. 28, SIAM, Philadelphia, 1999.}

\bibitem{GenPas}\foreignlanguage{english}{A. Gentile, A. Passarelli
di Napoli, \textit{Higher regularity for weak solutions to degenerate
parabolic problems}, Calc. Var. \textbf{62}, 225 (2023).}

\bibitem{Giann}\foreignlanguage{english}{F. Giannetti, A. Passarelli
di Napoli, C. Scheven, \textit{Higher differentiability of solutions
of parabolic systems with discontinuous coefficients}, J. Lond. Math.
Soc. (2) \textbf{94}, no. 1, 1-20 (2016).}

\bibitem{Giu}E. Giusti, \textit{Direct Methods in the Calculus of
Variations}, World Scientific Publishing Co., 2003.

\bibitem{Har}D.D. Haroske, \textit{Envelopes and Sharp Embeddings
of Function Spaces}, Chapman \& Hall CRC, 2006.

\bibitem{Irv}\foreignlanguage{english}{C. Irving, L. Koch, \textit{Boundary
regularity results for minimisers of convex functionals with $(p,q)$-growth},
Adv. Nonlinear Anal. \textbf{12}, no. 1 (2023).}

\bibitem{Isernia}T. Isernia, \textit{BMO regularity for asymptotic
parabolic systems with linear growth}, Differential Integral Equations
\textbf{28}, no. 11-12, 1173-1196 (2015).

\bibitem{Jones} P.W. Jones, \textit{Quasiconformal mappings and extendability
of functions in Sobolev spaces}, Acta Math. \textbf{147}, no. 1-2,
71-88 (1981).


\bibitem{Kuusi}\foreignlanguage{british}{T. Kuusi, G. Mingione, \textit{New
perturbation methods for nonlinear parabolic problems}, J. Math. Pures
Appl. (9) \textbf{98}, no. 4, 390-427 (2012).}


\bibitem{Leoni}G. Leoni, \textit{A First Course in Sobolev Spaces},
Grad. Stud. Math. 105, Amer. Math. Soc., Providence, 2009.

\bibitem{Lind1} P. Lindqvist,\textit{ On the time derivative in a
quasilinear equation}, Skr. K. Nor. Vidensk. Selsk.,\textbf{ 2}, 1-7
(2008).

\bibitem{Lind2} P. Lindqvist,\textit{ On the time derivative in an
obstacle problem}, Rev. Mat. Iberoam. \textbf{28}, no. 2, 577-590
(2012).

\bibitem{Lind3} P. Lindqvist,\textit{ The time derivative in a singular
parabolic equation}, Differential Integral Equations \textbf{30},
no. 9-10, 795-808 (2017).

\bibitem{LindPL} P. Lindqvist,\textit{ Notes on the Stationary $p$-Laplace Equation}, Springer Briefs Math., Springer, Cham, 2019.

\bibitem{Maz'ya}V. Maz'ya, \textit{Sobolev Spaces with Applications
to Elliptic Partial Differential Equations}, Grundlehren Math. Wiss.
342, Springer, Heidelberg, 2011.

\bibitem{Ru}\foreignlanguage{english}{S. Russo, \textit{On widely
degenerate $p$-Laplace equations with symmetric data}, Rev. Mat.
Complut. (2025). DOI: \url{https://doi.org/10.1007/s13163-025-00524-w}.}

\bibitem{Sche}\foreignlanguage{english}{C. Scheven, \textit{Regularity
for subquadratic parabolic systems: Higher integrability and dimension
estimates}, Proc. Roy. Soc. Edinburgh Sect. A \textbf{140}, no. 6,
1269-1308 (2010).}

\bibitem{Strunk}M. Strunk, \textit{Gradient regularity for widely
degenerate parabolic equations}, preprint (2025), \href{https://arxiv.org/abs/2510.07999}{arXiv: 2510.07999}. 

\bibitem{Tar}L. Tartar, \textit{An Introduction to Sobolev Spaces
and Interpolation Spaces}, Lecture Notes of the Unione Mat. Italiana,
vol. 3, Springer-Verlag, Berlin Heidelberg, 2007.

\bibitem{Tri0}\foreignlanguage{english}{H. Triebel, \textit{Spaces
of Besov-Hardy-Sobolev type}, Teubner-Texte Math. 15, Teubner, Leipzig,
1978.}

\bibitem{Tri}H. Triebel, \textit{Theory of Function Spaces}, Monogr.
Math. 78, Birkhäuser, Basel, 1983.

\bibitem{Tri2}H. Triebel, \textit{Theory of Function Spaces II},
Monogr. Math. 84, Birkhäuser, Basel, 1992. 

\bibitem{Uhl}K. Uhlenbeck, \textit{Regularity for a class of non-linear
elliptic systems}, Acta Math.\textbf{ 138}, no. 3-4, 219-240 (1977).

\bibitem{Zhou}Y. Zhou, \textit{Fractional Sobolev extension and imbedding},
Trans. Amer. Math. Soc. \textbf{367}, no. 2, 959-979 (2015).


\end{thebibliography}

\end{document}